\documentclass{amsart}
\usepackage[T1]{fontenc}
\usepackage[latin1]{inputenc}
\usepackage[french]{babel}
\usepackage{amssymb, amsmath, color, textcomp}
\usepackage[hypertex]{hyperref}
\renewcommand{\marginpar}[1]{}

\theoremstyle{plain}
\newtheorem{theorem}{Théorème}[section]
\newtheorem{corollary}[theorem]{Corollaire}
\newtheorem{prop}[theorem]{Proposition}
\newtheorem{lemma}[theorem]{Lemme}
\newtheorem*{thm}{Theorem}

\theoremstyle{definition}
\newtheorem{remark}[theorem]{Remarque}

\newtheorem{definition}[theorem]{D\'efinition}
\newtheorem{example}[theorem]{Exemple}

\def\satz{\begin{theorem}}
\def\esatz{\end{theorem}}
\def\satzli{\begin{prop}}
\def\esatzli{\end{prop}}
\def\bem{\begin{remark}}
\def\ebem{\end{remark}}
\def\kor{\begin{corollary}}
\def\ekor{\end{corollary}}
\def\lmm{\begin{lemma}}
\def\elmm{\end{lemma}}
\def\defn{\begin{definition}}
\def\edefn{\end{definition}}
\def\bsp{\begin{example}}
\def\ebsp{\end{example}}
\def\bew{\begin{proof}}
\def\ebew{\end{proof}}
\def\sstar{{}^\star_\star}

\newcommand{\Q}{\mathbb{Q}}
\newcommand{\cb}{\operatorname{Cb}}
\newcommand{\K}{\mathcal{K}}
\newcommand{\F}{\mathbb{F}}
\newcommand{\sscl}[1]{\langle #1 \rangle}

\newcommand{\tp}{\operatorname{tp}}
\newcommand{\lstp}{\operatorname{lstp}}
\newcommand{\RM}{\operatorname{RM}}
\renewcommand{\Im}{\operatorname{Im}}

\def\LL{\mathcal L}
\def\dcl{\mathrm{dcl}}
\def\acl{\mathrm{acl}}
\def\bdd{\mathrm{bdd}}
\def\p{\varphi}

\def\M{\mathfrak{M}}

\def\diag{\mathrm{diag}}
\def\ldim{\mathrm{dimlin}}
\def\NN{\mathcal N}

\def\Ind#1#2{#1\setbox0=\hbox{$#1x$}\kern\wd0\hbox to 0pt{\hss$#1\mid$\hss}
\lower.9\ht0\hbox to 0pt{\hss$#1\smile$\hss}\kern\wd0}
\def\Notind#1#2{#1\setbox0=\hbox{$#1x$}\kern\wd0\hbox to 0pt{\mathchardef\nn="0236\hss$#1\nn$\kern1.4\wd0\hss}\hbox to 0pt{\hss$#1\mid$\hss}\lower.9\ht0
\hbox to 0pt{\hss$#1\smile$\hss}\kern\wd0}
\def\ind{\mathop{\mathpalette\Ind{}}}
\def\nind{\mathop{\mathpalette\Notind{}}}
\def\indi#1{\mathop{\ \ \hbox to 0pt{\hss$\mid^{\hbox to 0pt{$\scriptstyle#1$\hss}}$\hss}
\lower4pt\hbox to 0pt{\hss$\smile$\hss}\ \ }}
\def\nindi#1{\mathop{\ \ \hbox to 0pt{\hss$\!\not{\mid}^{\hbox to 0pt{$\scriptstyle\,#1$\hss}}$\hss}
\lower4pt\hbox to 0pt{\hss$\smile$\hss}\ \ }}

\begin{document}

\title{Géométries relatives}
\date{\today}

\author{Thomas Blossier, Amador Martin-Pizarro et Frank O. Wagner}
\address{Universit\'e de Lyon; CNRS; Universit\'e Lyon 1; Institut Camille Jordan UMR5208, 43 boulevard du 11 novembre 1918, F--69622 Villeurbanne Cedex, France }
\email{blossier@math.univ-lyon1.fr}
\email{pizarro@math.univ-lyon1.fr}
\email{wagner@math.univ-lyon1.fr}
\thanks{Recherche soutenue par le réseau européen MRTN-CT-2004-512234 Modnet, le projet ANR-09-BLAN-0047 Modig, ainsi que l'Institut universitaire de France (F. Wagner).}
\keywords{Model Theory, Amalgamation methods, CM-triviality, Geometry, Group, Definability}

\subjclass{03C45}

\begin{abstract}  Une analyse des propriétés géométriques d'une structure relatives à un réduit est entamée. En particulier la définissabilité des groupes et des corps dans ce cadre est étudiée. Dans le cas relativement monobasé, tout groupe définissable est isogène à un sous-groupe d'un produit de groupes définissables dans les réduits. Dans le cas relativement CM-trivial, cas qui englobe certains amalgames de Hrushovski (la fusion de deux théories fortement minimales, les expansions d'un corps par un prédicat), tout groupe définissable s'envoie par un homomorphisme à noyau central dans un produit de groupes définissables dans les réduits.
\end{abstract}
\maketitle

\section*{English Summary} In this paper, we shall study type-definable groups in a simple theory with respect to one or several stable reducts. While the original motivation came from the analysis of definable groups  in structures obtained by Hrushovski's amalgamation method, the notions introduced are in fact more general, and in particular can be applied to certain expansions of algebraically closed fields by operators. We prove the following (Theorem \ref{groupe-homo}): 
\begin{thm} Let $T$ be simple and $T_0$ be a stable reduct of $T$. If $G$ is a type-definable group in $T$, there are a $*$-interpretable group $H$ in $T_0$ and a definable homomorphism $\phi:G^0\to H$ such that for independent generic elements $g$, $g'$ of $G$ we can name a set $D$ independent of $g$, $g'$ with
$$\acl(g),\acl(g')\indi0_{\acl(\phi(gg'))\cap\acl_0(\acl(\phi(g)),\acl(\phi(g')))}\acl(gg'),$$
where $\indi0$ and $\acl_0$ denote independence and algebraic closure respectively in the reduct $T_0$ over $D$.
\end{thm}
Thus $\phi$ captures all possible $0$-dependence relations resulting from the group operation. In order to ensure that the morphism is non-trivial, we need to add further hypotheses on the geometric complexity of $T$ over $T_0$: {\em Relative one-basedness} and {\em relative CM-triviality}. The kernel of $\phi$ is finite in the first case (Theorem \ref{groupe-monobase}) and virtually central in the latter (Theorem \ref{T:homo}).

Examples for relatively one-based theories are differentially closed fields in characteristic $0$ and fields with a generic automorphism. We therefore recover, up to isogeny, the characterisation of definable groups in these structures from \cite{KP02,Pi97}. Finally, we show that all coloured fields and fusions constructed by Hrushovski amalgamation are relatively CM-trivial (Theorem \ref{T:amalgames}).

\section{Introduction}

Un théorème de Pillay \cite{Pi97} affirme que tout groupe différentiellement constructible se plonge dans un groupe algébrique. Kowalski et Pillay \cite{KP02} ont obtenu un résultat analogue pour les groupes constructibles connexes sur un corps aux différences, modulo un noyau fini \cite{ChHr99}. Ces deux exemples nous montrent que des groupes {\em définissables} dans des structures {\em enrichies} peuvent \^etre analysés à partir de la structure de base. Dans les deux cas précédents, la démonstration consiste à passer de relations dans les corps enrichis à des relations purement algébriques dans une certaine configuration géométrique, dite {\em configuration de groupe}. \`A partir d'une telle configuration on récupère un groupe  \cite{eHPhD} qui sera ici algébrique. 

Il s'avère que des propriétés d'un groupe définissable dans une structure dépendent de l'existence de certains types de configurations géométriques. Un pseudo-plan généralise la notion d'incidence entre points et droites existant dans un plan euclidien. Une structure interprétant un pseudo-plan (complet) est {\em 1-ample} et cette configuration peut \^etre généralisée pour chaque $n$~: une géométrie {\em $n$-ample} correspond à l'existence d'un pseudo-espace en dimension $n$, ce qui nous donne toute une hiérarchie de complexités géométriques (hiérarchie conjecturée stricte, voir Evans \cite{Ev03}). Notons qu'un corps algébriquement clos est $n$-ample pour tout $n$.

Une structure est {\em monobasée} si elle n'est pas $1$-ample. D'après Hrushovski et Pillay \cite{HrPi} un groupe stable monobasé est abélien-par-fini et possède des propriétés de rigidité remarquables~: tout sous-ensemble définissable est une combinaison booléenne de translatés de sous-groupes définissables (en fait de sous-groupes définissables sur la clôture algébrique du vide, ce qui entraîne qu'il n'y a qu'un nombre borné de sous-groupes définissables). Cette propriété est à la base de la  démonstration modèle-théorique donnée par Hrushovski \cite{Hr96} de la conjecture de Mordell-Lang, qui réside dans la caractérisation des géométries associées aux types minimaux différentiels de variétés semi-abéliennes.

Dans le cas des structures non-2-amples (aussi appelées {\em CM-triviales} par Hrushovski \cite{Hr93}), Pillay \cite{Pi95} a montré que les groupes de rang de Morley fini sont nilpotents-par-fini~; ce résultat a été généralisé à certains groupes stables  \cite{Wa98}.

Dans cet article, nous définissons des notions de $1$-ampleur et $2$-ampleur relatives à un réduit par rapport à un opérateur de clôture. En imposant des propriétés supplémentaires à cet opérateur de clôture, nous réussissons à caractériser les groupes type-définissables dans de telles structures par rapport à ceux du réduit. On montre dans le cas relativement monobasé 
que tout groupe type-définissable a un sous-groupe d'indice borné qui se plonge définissablement modulo un noyau fini dans un groupe interprétable dans le réduit (Théorème \ref{groupe-monobase}). 
Les corps différentiellement clos et les corps munis d'un automorphisme générique sont des exemples de structures relativement monobasées (c.à.d.\ non $1$-amples) au-dessus du pur corps algébriquement clos sous-jacent. On retrouve ainsi(Corollaire \ref{KoPi-Pi}) les théorèmes de Kowalski-Pillay et de Pillay (quoique pour un groupe différentiel on obtienne seulement une monogénie, c'est-à-dire un plongement d'un sous-groupe d'indice fini modulo un noyau fini). Dans le cas relativement CM-trivial, à partir d'un groupe type-définissable
on construit un homomorphisme définissable, dont le noyau est virtuellement central, d'un sous-groupe d'indice borné dans un groupe interprétable dans le réduit (Théorème \ref{T:homo}). En particulier tout groupe simple type-défi\-nis\-sable se plonge dans un groupe du réduit (Corollaire \ref{cm-simple}), ce qui permet de montrer que tout corps type-définissable est définissablement isomorphe à un sous-corps d'un corps définissable dans le réduit (Corollaire \ref{C:corps}).

En utilisant la méthode d'amalgamation introduite par Hrushovski, diverses structures ayant des géométries exotiques furent construites: notamment des corps munis d'un prédicat pour un sous-ensemble (appelés corps {\em colorés} par Poizat) \cite{Po99, Po01, BaHo, BMPZ05,BHPW06}, ainsi que la fusion de deux théories fortement minimales \cite{Hr92, BMPZ06}, ou de rang de Morley fini et rang et degré de Morley définissables \cite{Z08}. Dans la dernière partie on vérifie que ces amalgames sont relativement CM-triviaux (c.à.d. non $2$-amples) au dessus des théories de base.
Notons que la CM-trivialité relative ne permet pas d'analyser tous les groupes définissables. Dans les corps colorés non-collapsés, on vérifie facilement que tout groupe définissable connexe est l'extension d'un sous-groupe coloré par un groupe algébrique. Ceci n'est plus vrai pour certains groupes interprétables, comme le quotient du corps entier par le sous-groupe coloré. Notons que dans le cas collapsé de rang fini, tout groupe interprétable devient définissable par élimination des imaginaires \cite{Wa01}.

Quant à la fusion, Hrushovski avait défini dans \cite{Hr93} une notion de {\em platitude} (absolue) qui empêche l'existence de groupes~; il affirme dans \cite{Hr92} que la fusion sur un réduit commun trivial est {\em plate} relativement aux théories de base, et que cela implique que tout groupe définissable est isomorphe à un produit direct, à centre fini près, de deux groupes définissables dans les deux théories de départ. Dans un travail en cours nous utilisons cette notion afin d'étudier les groupes abéliens.

Bien que ce papier demande une certaine connaissance d'outils modèle-théoriques (voir \cite{Po85}, \cite{Pi96} pour la stabilité et \cite{wa00} pour la simplicité), une connaissance approfondie de la méthode d'amalgamation n'est pas nécessaire pour sa lecture (la partie \ref{S:Amal} est indépendante et contient un rappel de cette méthode). 

Nous remercions A. Pillay pour ses commentaires qui nous ont permis de démar\-rer ce travail, et le rapporteur pour sa lecture critique et détaillée.

\section{Prélude}

Dans cet article nous allons considérer une théorie complète simple $T$ dans un langage $\LL$ avec un réduit stable $T_0$ à un sous-langage $\LL_0$, ou encore avec une famille de réduits stables $(T_i:i<n)$ à des sous-langages $\LL_i$. Les notions modèle-théoriques comme la clôture définissable $\dcl$, la clôture algébrique $\acl$, la clôture bornée $\bdd$, les types $\tp$, les bases canoniques $\cb$ ou l'indépendance $\ind$ s'entendent au sens de $T$~; si on les prend au sens de $T_i$ on l'indiquera par l'indice $i$: $\dcl_i$, $\acl_i$, $\tp_i$, $\cb_i$, $\indi i$.

Rappelons que si $E$ est une relation d'équivalence définissable sur un uplet fini $a$, alors la classe $a_E$ de $a$ modulo $E$ est un {\em imaginaire}; si $E$ n'est que type-définissable et $a$ un uplet possiblement infini (appel\'e {\em $*$-uplet}), alors $a_E$ est un {\em hyperimaginaire}. En particulier les $*$-uplets eux-m\^emes sont des hyperimaginaires. Les imaginaires possèdent les mêmes propriétés modèle-théoriques que les éléments et uplets ordinaires (appel\'es {\em r\'eels}); ceci vaut aussi pour les hyperimaginaires, sauf que l'\'egalit\'e de deux hyperimaginaires est un type partiel plut\^ot qu'une formule. Un ensemble est {\em d\'efinissable} s'il est donné par une formule, et {\em type-d\'efinissable} s'il est donn\'e par un type partiel d'arit\'e finie; il est {\em (type-)interprétable} si la formule ou le type partiel portent sur des variables imaginaires; un  préfixe $*$- indiquera que le nombre de variables (et de formules) est infini. Ainsi, si $G$ est un groupe type-définissable (le domaine et le graphe de la multiplication le sont), on voit facilement que le graphe de la multiplication est en fait donné par une seule formule; ceci n'est plus vrai pour un groupe $*$-définissable. Enfin, si les réduits éliminent complètement les imaginaires, alors interprétable équivaut à définissable; nous avons maintenu la distinction afin d'indiquer quand on aborde l'imaginaire.

Remarquons que si $E$ est une relation d'équivalence type-définissable dans $T$ mais pas dans $T_i$, alors une classe modulo $E$ n'a aucun sens dans $T_i$~; les hyperimaginaires de $T$ n'existent pas nécessairement dans les réduits. On travaillera donc sauf mentions contraires uniquement avec des éléments réels~; en particulier les clôtures algébrique et définissable sont restreintes aux réels (même si on prend la clôture d'un ensemble d'hyperimaginaires). Par contre, la clôture bornée s'entend toujours au sens hyperimaginaire. 

Afin d'avoir néanmoins un certain contrôle sur les imaginaires, nous supposerons par la suite que les $T_i$ éliminent géométriquement les imaginaires pour tout $i<n$, c'est-à-dire tout $T_i$-imaginaire est $T_i$-interalgébrique avec un uplet réel. Ceci est satisfait par exemple si les $T_i$ sont fortement minimales avec $\acl_i(\emptyset)$ infini. Notons en particulier que si nous travaillons avec un seul réduit $T_0$, l'élimination géométrique de imaginaires est obtenu gratuitement en ajoutant les $0$-imaginaires au langage. 

Les énoncés principaux de cet article, notamment les théorèmes \ref{groupe-homo}, \ref{groupe-monobase} et \ref{T:homo}, restent valables si le groupe de départ est $*$-définissable. Dans ce cas les groupes obtenus dans les réduits sont $*$-interprétables.

Plus généralement, on pourrait partir de réduits qui ne sont que simples avec élimination géométrique des hyperimaginaires et utiliser le théorème de configuration de groupe pour les théories simples \cite{BYTW}. Nous ne voyons pas d'obstacle majeur à une extension de nos résultats à ce contexte, mais la non-unicité des extensions non-déviantes et le travail avec les presque-hyperimaginaires nécessiteraient des arguments supplémentaires que nous n'avons pas tenté de faire en détail.

Le lemme suivant, moins évident qu'il ne le semble à première vue, sera utilisé fréquemment dans cet article.

\lmm\label{L:indep}
Si $B$ est algébriquement clos (au sens de $T$) et $a\ind_B c$, alors $a\indi 0_B c$.
\elmm
\bew Comme la théorie $T_0$ est stable, il suffit de voir que $\tp_0(a/B c)$ ne divise pas sur $B$. Fixons $\varphi(x,y)$ une $\LL_0$-formule à paramètres sur $B$ satisfaite par $(a,c)$.  Considérons une suite de Morley $(c_j:j <\omega\cdot 2)$ de $\tp(c/B)$. Alors $\bigwedge_i\p(x,c_i)$ est consistant, comme $a$ et $c$ sont indépendants sur $B$. Pour montrer que $\p(x,c)$ ne $0$-divise pas, il suffit  de montrer que $(c_j:j<\omega\cdot 2)$ est une $0$-suite de Morley. Or, la $0$-indiscernabilité ne posant aucun problème, il s'agit de montrer qu'elle est $0$-indépendante sur $B$.

Remarquons que $(c_j:\omega\le j<\omega\cdot2)$ est $0$-indépendante sur $(B,c_j:j<\omega)$ par stabilité, et forme donc une $0$-suite de Morley de $\tp_0(c_\omega/B,c_j:j<\omega)$.
Puisque la base canonique d'un type est algébrique sur une suite de Morley de réalisations, on a
$$\cb_0(c_\omega/B,c_j:j<\omega)\subseteq\acl_0^{eq}(B,c_j: j<\omega)\cap\acl_0^{eq}(c_j:\omega\le j<\omega\cdot 2).$$
Or, $T_0$ élimine géométriquement les imaginaires, $B$ est algébriquement clos, et
$$(c_j: j<\omega)\ind_B (c_j:\omega\le j<\omega\cdot 2).$$
Ceci implique que $\cb_0(c_\omega/B,c_j:j<\omega)\subseteq\acl_0^{eq}(B)$, et donc
$$c_\omega\indi 0_B(c_j:j<\omega).$$
Par indiscernabilité, la suite $(c_j:j<\omega\cdot 2)$ est bien $0$-indépendante sur $B$.\ebew

\section{Groupes et réduits}

Cette partie est consacrée à la démonstration du théorème suivant, qui établit l'existence d'un homomorphisme d'un groupe $T$-définissable vers un groupe interprétable dans un réduit $T_0$. Cet homomorphisme peut s'avérer trivial sans conditions supplémentaires, comme les non-ampleurs relatives qui seront discutées dans les parties suivantes. 

\satz\label{groupe-homo} Soit $T$ une théorie simple avec un réduit $T_0$ stable avec élimination géométrique des imaginaires.
Tout groupe $G$ type-définissable sur $\emptyset$ dans $T$ a un sous-groupe normal $N$ type-définissable sur $\emptyset$ et un sous-groupe type-définissable $G_N$ d'indice borné dans $G$ tel que $G_N/N$ se plonge d\'efinissablement (\`a l'aide de de paramètres $A$ éventuels) dans un groupe $H$ $*$-interprétable dans $T_0$. De plus, pour tout $B\supseteq A$ et tous génériques indépendants $g_0,g_1$ de $G$ sur $B\supseteq A$ il existe $D\supseteq B$ indépendant de $g_0,g_1$ sur $B$ avec
$$\acl(g_0,D),\acl(g_1,D)\indi0_{\acl(g_0g_1N,D)\cap\acl_0(\acl(g_0N,D),\acl(g_1N,D))}\acl(g_0g_1,D).$$
\esatz

\bem\label{homo-prec} Rappelons que par définition $\acl(gN,D)$ n'est constitué que d'éléments réels. Mais on choisira $A$ tel que $gN$ soit interalgébrique sur $A$ avec un $*$-uplet réel. En particulier $gN\in\bdd(\acl(gN,D))$.

Notons aussi que $G_N/N$ est un hyperimaginaire: Les éléments du quotient sont des classes modulo la relation d'équivalence type-définissable $x^{-1}y\in N$. Un plongement définissable dans $H$ sera donc donn\'e par un homomorphisme $*$-interprétable de $G_N$ dans $H$ dont le noyau est \'egal \`a $N$. (Le graphe d'un tel homomorphisme est l'image d'un ensemble $*$-définissable de $G_N\times X$ par la projection de $G_N\times X$ sur $G_N\times H$, où $H=X/E$ pour un ensemble $*$-définissable $X$ et une relation d'équivalence $*$-définissable $E$.) Un tel plongement devrait être dit  {\em hyperdéfinissable}, mais nous ne souhaitons pas alourdir la nomenclature.\ebem

\bew Soit $G$ un groupe $T$-type-définissable sur $\emptyset$. Rappelons que deux sous-groupes $H$ et $K$ de $G$ sont dits \emph{commensurables} si les indices de $H\cap K$ dans $H$ et $K$ sont bornés. On appellera un sous-groupe type-définissable (avec paramètres) $N\le G$ {\em nucléaire} s'il est normal dans un sous-groupe type-définissable $G_N$ d'indice borné dans $G$ tel que $G_N/N$ se plonge définissablement dans un groupe $H_N$ $*$-interprétable dans $T_0$.

Pour cette démonstration nous vérifions deux résultats préliminaires indépendants (Lemmes \ref{noyauminimal} et \ref{reduitgroupe}) que nous combinons ensuite pour conclure.

\lmm\label{noyauminimal} Il existe un sous-groupe nucléaire $N$ de $G$ type-définissable sur $\emptyset$ et minimal à commensurabilité près parmi les sous-groupes nucléaires de $G$.\elmm

\bew Soit $\NN$ la famille des sous-groupes nucléaires de $G$~; notons qu'elle est non-vide puisqu'elle contient $G$. Elle est évidemment invariante par $\emptyset$-automor\-phisme, et par automorphisme définissable de $G$ (et en particulier par conjugaison). De plus, elle est stable par intersection bornée : si $\NN_0$ est un sous-ensemble borné de $\NN$ alors  l'intersection $\bigcap_{N\in\NN_0}N$  est type-définissable, normal dans $\bigcap_{N\in\NN_0}G_N$, et $\bigcap_{N\in\NN_0}G_N/\bigcap_{N\in\NN_0}N$ se plonge définissablement dans le groupe $\prod_{N\in\NN_0}H_N$, qui est $*$-interprétable dans $T_0$. 

Si la théorie $T$ est stable, l'intersection $\bigcap_{N\in\NN}N$ est \'egale \`a une sous-intersection de taille au plus $|T|$ et donc type-définissable (voir \cite[Corollaire 1.0.8]{Wa97}); elle est $\emptyset$-invariante puisque $\NN$ l'est. D'apr\`es le paragraphe précédent elle appartient à $\NN$, ce qui permet de conclure.

Plus généralement d'après \cite[Theorem 4.5.13]{wa00} appliqué à $\NN$ il existe un sous-groupe $N$ normal dans $G$ et type-définissable sur $\emptyset$ contenant un sous-groupe d'indice borné  $N_1 \in\NN$ et tel que chaque intersection d'un $N_2\in\NN$ avec $N$ est d'indice borné dans $N$. Vérifions pour terminer que $N$ est nucléaire. Posons $G_N = N G_{N_1}$. Alors, le quotient $G_N/N$ est isomorphe au quotient $(G_{N_1}/N_1)/((N \cap G_{N_1})/N_1)$. Soit $H_{N_1}$ un groupe $*$-interprétable dans $T_0$ et $\phi$ un plongement définissable de $G_{N_1}/N_1$ dans $H_{N_1}$. Comme $(N \cap G_{N_1})/N_1$ est borné, son image $I$ dans $H_{N_1}$ l'est aussi. Donc $I$ est une limite projective de groupes finis, et en particulier $*$-interprétable dans $T_0$. Alors $N_{H_{N_1}}(I)$ est aussi $*$-interprétable dans $T_0$. Donc $H_N=N_{H_{N_1}}(I)/I$ est $*$-interprétable dans $T_0$, et $\phi$ induit un plongement définissable de $G_N/N$ dans $H_N$.
\ebew

Soit $\phi$ le plongement de $G_N/N$ dans un groupe $H$ $*$-interprétable dans $T_0$, et $A$ les paramètres nécessaires pour définir $G_N$, $H$ et $\phi$. On complète $A$ avec un système de représentants pour $G/G_N$. Ainsi pour tout $g\in G$ il y a $g'\in A$ avec $g'g\in G_N$, et le translaté $gN$ est interdéfinissable sur $A$ avec $g'gN$ et donc par élimination géométrique des imaginaires dans $T_0$, interalgébrique avec $\acl(\phi(g'g),A)$.

Soit $B$ un ensemble de paramètres contenant $A$. Prenons trois éléments génériques $g_0$, $g_1$ et $g_3$ de $G$ indépendants sur $B$. Il nous faut maintenant trouver $D$ vérifiant la condition sur la $0$-indépendance.

Par l'interdéfinissabilité précédente, on peut supposer pour la suite que $g_0$, $g_1$ et $g_3$ sont dans $G_N$. Afin d'alléger les notations on ajoute $B$ au langage. Posons $g_2=g_0 g_1$, $g_4=g_2  g_3$ et $g_5=g_1 g_3 = g_0^{-1}g_2g_3=g_0^{-1}g_4$. Alors le $6$-uplet $(g_0,g_1,g_2,g_3,g_4,g_5)$ 
\begin{center}
\begin{picture}(50,70)
\put(0,0){\line(1,3){20}}
\put(20,60){\line(4,-1){68}}
\put(10,30){\line(6,1){78}}
\put(0,0){\line(1,1){52}}
\put(-12,-2){$g_2$}
\put(-3,28){$g_1$}
\put(8,58){$g_0$}
\put(50,55){$g_4$}
\put(33,25){$g_3$}
\put(90,39){$g_5$}
\end{picture}
\end{center}
forme une configuration de groupe, c'est-à-dire les éléments de chaque paire ainsi que les éléments de chaque triplet non-colinéaire sont indépendants, mais sur toute droite chaque point est algébrique sur les deux autres.

\lmm\label{reduitgroupe} Il existe $D$, un ensemble dénombrable de génériques de $G$ indépendants sur $g_0,g_1,g_3$, tel que pour tout triplet $(g_i,g_j,g_k)$ colinéaire avec $0\le i,j,k\le 5$, l'intersection
$$\alpha_i=\alpha_i(j,k)=\acl(g_i,D)\cap\acl_0(\acl(g_j,D),\acl(g_k,D))$$
ne dépend pas du choix de $g_j,g_k$. De plus
$$\acl(g_j,D),\acl(g_k,D)\indi0_{\alpha_i}\acl(g_i,D).$$
En particulier, $(\alpha_0,\alpha_1,\alpha_2,\alpha_3,\alpha_4,\alpha_5)$ forme une $0$-configuration de groupe sur $\acl(D)$.\elmm

\bew 
Considérons $D_0$ un ensemble dénombrable d'éléments génériques de $G$ indépendants sur $g_0,g_1,g_3$. Notons que $g_0$ est générique sur $g_1, D_0$ et $g_3$ l'est sur $g_0,g_1,D_0$. Si $I$ est une suite de Morley de $\tp(g_0/\acl(g_1,D_0))$ indépendante de $g_0,g_3$, ou une suite de Morley  de $\tp(g_3/\acl(g_0,g_1,D_0))$ indépendante de $g_3$, alors les ensembles $D_0\cup g_0 I$ et $D_0\cup g_0^{-1} I$ sont toujours des ensembles d'éléments génériques indépendants sur $g_0,g_1,g_3$. En itérant sur les triplets de points non colinéaires du diagramme il est ainsi possible d'étendre $D_0$ en un ensemble dénombrable $D_1$ d'éléments génériques indépendants sur $g_0,g_1,g_3$, tel que pour tout couple de points $(g_i,g_j)$ sur une même droite, et  tout autre point $g_\ell$ non colinéaire,
\begin{enumerate}
\item les translatés à gauche $g_i  D_1$ et $g_i^{-1}  D_1$ contiennent chacun une suite de Morley de  $\tp(g_i/\acl(g_j, D_0))$, et
\item $D_1$ contient  une suite de Morley de $\tp(g_\ell/\acl(g_i,g_j, D_0))$.\end{enumerate}

En itérant $\omega$ fois cette construction on obtient un ensemble dénombrable $D$ d'éléments génériques indépendants sur $g_0,g_1,g_3$, tel que pour tout uplet fini $\bar d$ de $D$, tout couple de points $(g_i,g_j)$ sur une même droite et  tout autre point $g_\ell$ non colinéaire,

\begin{enumerate}
\item $g_i  D$ et $g_i^{-1}  D$ contiennent tous deux des suites de Morley de  $\tp(g_i/\acl(g_j, \bar d))$, et
\item $D$ contient  une suite de Morley  de $\tp(g_\ell/\acl(g_i,g_j,\bar d))$.\end{enumerate}

Pour trois points $g_i,g_j,g_k$ sur une même droite et un uplet fini  $\bar d$ de $D$, on pose 
$$\alpha(i;j,k,\bar d)=\acl_0(\cb_0(\acl(g_j,\bar d),\acl(g_k,\bar d)/\acl(g_i,\bar d))),$$
un $*$-uplet réel $0$-interalgébrique avec $\cb_0(\acl(g_j,\bar d),\acl(g_k,\bar d)/\acl(g_i,\bar d))$ par élimination géométrique des imaginaires dans $T_0$.

Afin d'alléger les notations, supposons que $g_ig_j=g_k$. (Dans les autres cas on remplace $g_i$, $g_j$ ou $g_k$ par son inverse, qui lui est $T$-interalgébrique.) Par construction, $g_j D$ contient une suite de Morley $I$ de $\tp(g_j/\acl(g_i,\bar d))$. Alors, la suite
$$(\acl(e,\bar d),\acl(g_i   e ,\bar d):e \in I)$$
est une suite de Morley de $\tp(\acl(g_j,\bar d),\acl(g_k,\bar d)/\acl(g_i,\bar d))$. D'après le Lemme \ref{L:indep} elle est aussi $0$-indépendante sur $\acl(g_i,\bar d)$, c'est-à-dire est une $0$-suite de Morley. Donc
$$\alpha(i;j,k,\bar d)\subseteq \acl_0(\acl(e,\bar d),\acl(g_i  e,\bar d):e\in I).$$
Mais pour tout $e\in I$, on a $e\in g_j D\subseteq\acl(g_j,D)$ et $g_i e\in g_i g_j D=g_k D\subseteq\acl(g_k,D)$. Donc 
$$\alpha(i;j,k, \bar d)\subseteq \acl_0(\acl(g_j,D),\acl(g_k,D))\cap\acl(g_i,D)=\alpha_i(j,k).$$
De plus, $g_i,g_j,g_k\ind_{\bar d}D$ pour tout uplet fini $\bar d$ dans $D$, ce qui implique d'après le Lemme \ref{L:indep} que
$$\acl(g_j,\bar d),\acl(g_k,\bar d)\indi0_{\acl(g_i,\bar d)}\acl(g_i,D)$$
et donc
\begin{align*}\cb_0(\acl(g_j,D),\acl(g_k,D)/\acl(g_i,D))&=\bigcup_{\bar d\in D}\cb_0(\acl(g_j,\bar d),\acl(g_k,\bar d)/\acl(g_i,D))\\
&=\bigcup_{\bar d\in D}\cb_0(\acl(g_j,\bar d),\acl(g_k,\bar d)/\acl(g_i,\bar d))\\ 
&\subseteq\bigcup_{\bar d\in D}\acl_0^{eq}(\alpha(i;j,k,\bar d))\subseteq\acl_0^{eq}(\alpha_i(j,k)).\end{align*}
Comme $\alpha_i(j,k)=\acl(g_i,D)\cap\acl_0(\acl(g_j,D),\acl(g_k,D))$, on en déduit que
\begin{equation}\label{eqn1}
\acl(g_j,D),\acl(g_k,D)\indi0_{\alpha_i(j,k)}\acl(g_i,D) 
\end{equation}
et $\alpha_i(j,k)=\bigcup_{\bar d\in D}\alpha(i;j,k,\bar d)$. 

Pour montrer que $\alpha_i(j,k)$ ne dépend pas de la droite choisie, considérons l'autre droite passant par $g_i$ donnée par les points $g_i$,  $g_m$, $g_n$. Soit $g_\ell$ le sixième point. 
On supposera $$g_j g_\ell=g_n\quad\mbox{et}\quad g_k g_\ell^{-1}=g_m,$$
les autres cas étant analogues.
Fixons un uplet fini $\bar d$ de $D$. Par construction, $D$ contient une suite de Morley $J$ de $\tp(g_\ell/\acl(g_i,g_j,\bar d))$. La suite
$$(\acl(g_j e,\bar d),\acl(g_k e^{-1},\bar d):e\in J)$$
est alors une suite de Morley de $\tp(\acl(g_n,\bar d),\acl(g_m,\bar d)/\acl(g_i,\bar d))$. Par le Lemme \ref{L:indep} elle est aussi $0$-indépendante sur $\acl(g_i,\bar d)$, et 
$$\alpha(i;m,n,\bar d)\subseteq\acl_0(\acl(g_j e,\bar d),\acl(g_k e^{-1},\bar d):e \in J).$$
Puisque $g_j e\in\acl(g_j,D)$ et $g_k e^{-1}\in\acl(g_k,D)$ pour tout $e \in J$, on en déduit que
$$\alpha(i;m,n,\bar d)\subseteq\acl(g_i,D)\cap\acl_0(\acl(g_j,D),\acl(g_k,D))=\alpha_i(j,k).$$
Donc $\alpha_i(m,n)\subseteq\alpha_i(j,k)$ et par symétrie on a l'égalité.

Montrons maintenant que  $\alpha_i$ est $0$-algébrique sur $\alpha_j, \alpha_k$. Par définition $\alpha_i$ est inclus dans $\acl_0(\acl((g_j,D),\acl(g_k,D))$. L'équation (\ref{eqn1}) nous donne l'indépendance
$$\acl(g_i,D),\acl(g_k,D)\indi0_{\alpha_j}\acl(g_j,D),$$
ce qui entraîne $\alpha_i\subset \acl_0(\alpha_j,\acl(g_k,D))$. 
Mais l'équation (\ref{eqn1}) donne aussi l'indépendance
$$\alpha_i,\alpha_j\indi0_{\alpha_k}\acl(g_k,D)$$
qui implique que $\alpha_i\subset \acl_0(\alpha_j,\alpha_k)$.
Enfin, dans la configuration $(\alpha_0,\alpha_1,\alpha_2,\alpha_3,\alpha_4,\alpha_5)$ l'indépendance des éléments des paires et des triplets non-colinéaires suit des indépendances de la configuration originelle~: Comme $(g_0,g_1,g_2,g_3,g_4,g_5)$ est indépendant de $D$, ça reste une configuration de groupe sur $\acl(D)$. Or, $\alpha_s\in\acl(g_s,D)$ pour $0\le s\le 5$. 
On obtient donc une $0$-configuration de groupe sur $\acl(D)$.
\ebew
\bem\label{R:DT_i} Pour chaque ensemble de paramètres $B \supseteq A$ et tous génériques indépendants $g_0$, $g_1$ de $G$ sur $B$, %thomas
l'ensemble $D\setminus B$ est constitué de génériques indépendants sur $B,g_0,g_1$ et sa construction 
ne dépend pas du réduit $T_0$. Si $T$ est stable et $G$ est connexe, par unicité du type générique on obtient $D$ en prenant une suite de Morley dénombrable du type générique de $G$ au-dessus de $g_0,g_1,g_3, B$.\ebem

Nous allons maintenant utiliser la $0$-configuration de groupe $(\alpha_0,\alpha_1,\alpha_2,\alpha_3,\alpha_4,\alpha_5)$ sur $\acl(D)$ obtenue afin de vérifier la condition de $0$-indépendance. D'après \cite{eHPhD} ou \cite[Theorem 5.4.5 et Remark 5.4.9]{Pi96} (où c'est fait pour les uplets finis, mais la démonstration se transfère {\em verbatim} au cas des $*$-uplets en remplaçant définissable par $*$-définissable) 
il existe un groupe $H'$ connexe $*$-interprétable dans $T_0$ sur des paramètres $C=\acl(C)\supseteq\acl(D)$ avec
$$C\indi0_{\acl(D)}\alpha_0,\alpha_1,\alpha_2,\alpha_3,\alpha_4,\alpha_5,$$
et trois génériques $0$-indépendants $h_0$, $h_1$ et $h_3$ de $H'$ avec $h_2=h_0h_1$, $h_4=h_0h_1h_3$ et $h_5=h_1h_3$, tels que $\alpha_j$ et $h_j$ sont $0$-interalgébriques sur $C$ pour chaque $j$ dans $\{0,\ldots,5\}$. Notons que la stabilité de $T_0$ implique que $H'$ est en fait une limite projective de groupes $T_0$-interprétables.

Par le Lemme \ref{L:indep} et la stationnarité de $\tp_0(C/\acl(D))$ on peut même supposer que $C\ind_D g_0,g_1,g_3$, ce qui implique que $h_0,h_1,h_3$ sont indépendants sur $C$ au sens de $T$.

Rappelons que le stabilisateur à droite de $\tp(g_0,h_0/C)$ dans $G\times H'$ est un sous-groupe $*$-interprétable $S$. L'ensemble 
$$E=\left\{(g,h)\in G\times H'\ :\ \parbox{20em}{$\exists\,(x,y)\models\tp(g_0,h_0/C)\ \big[(x,y)\ind_C(g,h)\ \land$\\
$(xg,yh)\models\tp(g_2,h_2/C)\land(xg,yh)\ind_C(g,h)\big]$}\right\}.$$
est générique dans un translaté à droite $*$-interprétable sur $C$ de $S$. Comme $E$ contient $(g_1,h_1)$ et $g_1$ est générique dans $G$ sur $C$, la projection de $S$ sur la première coordonnée est un sous-groupe $G_S$ d'indice borné dans $G$. Par contre la projection sur la deuxième coordonnée n'est pas nécessairement générique~: le type $\tp(h_1/C)$ n'est pas forcément générique dans $H'$ bien que $\tp_0(h_1/C)$ le soit au sens de $T_0$. Comme $h_0$ est algébrique sur $C,g_0$, le sous-groupe $H_0= \{h\in H':(1,h)\in S\}$ est une limite projective (de cardinalité bornée) de sous-groupes finis, donc $*$-définissable dans $T_0$, et le quotient $\Im(S)=N_{H'}(H_0)/H_0$ est $*$-interprétable dans $T_0$. Le noyau $\ker(S)=\{g\in G_S:(g,1)\in S\}$ est un sous-groupe distingué de $G_S$, et $S$ induit un plongement de $G_S/\ker(S)$ dans $\Im(S)$. On conclut ainsi que $\ker(S)$ est nucléaire.

Nous allons maintenant comparer le morphisme induit par $S$ avec le morphisme $\phi$ déterminé par le sous-groupe nucléaire minimal $N$.
Par minimalité, l'intersection $\ker(S)\cap N$ est d'indice borné dans $N$ et donc pour tout $g$ de $G$ la classe $g \ker(S)$ est dans $\bdd(gN,C)$. Le choix de $A$ assure que pour tout $g$ de $G_N$, l'image $\phi(g)$ est algébrique sur $g,A$ et donc sur $g,D$.

En particulier, pour $\{i,j,k\}=\{0,1,2\}$, on a 
$$\acl_0(\phi(g_i),A)\subseteq\acl(g_i,D)\cap\acl_0(\acl(g_j,D),\acl(g_k,D))=\alpha_i.$$
Donc
\begin{align*}\bdd(g_iN,A)&=\bdd(\phi(g_i),A)\subseteq\bdd(\alpha_i)\\
&\subseteq\bdd(h_i,C)=\bdd(g_i\ker(S),C)\subseteq\bdd(g_iN,C).\end{align*}
On en déduit que $g_iN$ et $\alpha_i$ sont interbornés sur $C$, et donc sur $D$ puisque $$g_iN,\alpha_i\ind_DC.$$
Comme $\alpha_i$ est un uplet réel, on a $\alpha_i \subset \acl(g_iN,D)$. 
Donc
\begin{align*} \alpha_2 &\subseteq\acl(\alpha_2)\cap \acl_0(\acl(\alpha_0),\acl(\alpha_1))\\ &\subseteq\acl(g_2N,D)\cap\acl_0(\acl(g_0N,D),\acl(g_1N,D)) \subseteq \alpha_2,\end{align*}
ce qui termine la démonstration du théorème.\ebew

Dans ce théorème le groupe $H$ peut être trivial, par exemple si on considère le réduit à l'égalité. Afin de contrôler le noyau de l'homomorphisme, nous allons introduire des conditions géométriques relatives.

\section{Théories relativement monobasées}

Le but originel de ce travail était l'étude des groupes définissables dans les amalgames. Dans ces constructions la clôture autosuffisante joue un rôle central (voir la partie \ref{S:Amal}). Nous supposerons donc pour la suite que la théorie $T$ de départ est munie d'un opérateur de clôture $\sscl.$ finitaire et invariant tel que $A\subseteq \sscl A\subseteq\acl(A)$ pour tout ensemble réel $A$.

\begin{definition}\label{defonebased}
La théorie $T$ est \emph{monobasée au dessus de $T_0$ pour $\sscl.$} si pour tous ensembles réels $A\subseteq B$ algébriquement clos et tout uplet réel $\bar c$, si 
$$\sscl{A\bar c} \indi0_{A} B,$$
alors la base canonique $\cb(\bar c/B)$ est bornée sur $A$. (Notons que cette base canonique est un hyperimaginaire de $T$).\newline
Plus généralement, $T$ est \emph{monobasée au dessus de $(T_i:i<n)$ pour $\sscl.$} si pour tout $A\subseteq B$ algébriquement clos et tout uplet $\bar c$, si
$$\sscl{A\bar c}\indi i_A B\mbox{ pour tout }i<n,$$
alors la base canonique $\cb(\bar c/B)$ est bornée sur $A$.\end{definition}

\begin{remark}
Toute th\'eorie est monobasée au dessus d'elle même pour l'opéra\-teur $\acl$. Pour le même opérateur de clôture, si une th\'eorie $T$ est monobasée au-dessus de son réduit au langage de l'égalité, alors elle est monobasée au sens classique ; la réciproque est vraie si $T$ élimine géométriquement les hyperimaginaires. Bien sur, si $T$ est stable, la base canonique $\cb(\bar b/B)$ est constituée d'éléments imaginaires de $T$ et elle est bornée sur $A$ si et seulement si elle est algébrique sur $A$.\end{remark}

\begin{definition}
La th\'eorie $T$ est \emph{$1$-ample au dessus de $T_0$ pour $\sscl.$} s'il existe des uplets réels $\bar a$, $\bar b$ et $\bar c$ tels que:
\begin{itemize}
\item $ \acl(\bar a,\bar b) \indi0_{\acl(\bar a)} \sscl{\acl(\bar a),\bar c}$.
\item $ \bar c\nind_{\bar a}\bar b$.
\end{itemize}
Plus généralement, $T$ est \emph{$1$-ample au dessus de $(T_i:i<n)$ pour $\sscl.$} s'il existe des uplets $\bar a$, $\bar b$ et $\bar c$ tels que:
\begin{itemize}
\item $\acl(\bar a,\bar b) \indi i_{\acl(\bar a)} \sscl{\acl(\bar a),\bar c}$ pour tout $i<n$.
\item $\bar c\nind_{\bar a}\bar b$.
\end{itemize}
\end{definition}
\bem
Une th\'eorie $T$ est monobasée au-dessus de $(T_i:i<n)$ pour $\sscl.$ si et seulement si elle n'est pas $1$-ample au dessus de $(T_i:i<n)$.\ebem
\bew On prend $A=\acl(\bar a)$ et $B=\acl(\bar a,\bar b)$.\ebew

Pour caractériser les groupes définissables dans des théories relativement monobasées, nous imposerons des conditions supplémentaires entre l'opérateur $\sscl.$ et les clôtures algébriques $\acl_i$~:
\begin{itemize}
\item[(\dag)] Si $A$ est algébriquement clos et $b\ind_Ac$, alors $\sscl{Abc}\subseteq\bigcap_{i<n}\acl_i(\sscl{Ab},\sscl{Ac})$.
\item[(\ddag)] Si $\bar a\in\bigcup_{i<n}\acl_i(A)$, alors $\sscl{\acl(\bar a),A}\subseteq\bigcap_{i<n}\acl_i(\acl(\bar a),\sscl{A})$.
\end{itemize}
\bem
Notons que la cl\^oture algébrique d'une théorie ne satisfait pas la condition $(\dag)$ au-dessus du réduit à l'égalité dès qu'un groupe non-trivial est définissable.
\ebem

\satzli\label{mb-para} Supposons que $\sscl.$ satisfasse $(\dag)$. Alors la $1$-ampleur relative est conservée par l'adjonction ou la suppression de paramètres.\esatzli
\bew C'est évident pour la suppression de paramètres. Pour l'adjonction, considérons un ensemble $D$ de paramètres tel que $T(D)$ soit monobasée au dessus de $(T_i:i<n)$. Soient $A\subseteq B$ algébriquement clos, et un uplet $\bar c$ tels que $\sscl{A\bar c}\indi i_AB$ pour tout $i<n$. On peux toujours supposer $B\bar c\ind_A D$. Soit $A'=\acl(AD)$ et $B'=\acl(BD)$. Alors
$\sscl{A\bar c}\ind_B B'$~; donc $\sscl{A\bar c}\indi i_B B'$ pour tout $i<n$, d'o\`u $\sscl{A\bar c}\indi i_A B'$ par transitivité. Donc $\sscl{A\bar c}\indi i_{A'}B'$. Puisque $\sscl.$ satisfait (\dag), l'indépendance $\bar c\ind_AA'$ implique $\sscl{A'\bar c}\subseteq\acl_i(\sscl{A\bar c},A')$ pour tout $i<n$, et donc $\sscl{A'\bar c}\indi i_{A'}B'$.

Comme $T(D)$ est relativement monobasée, $\cb(\bar c/B')$ est bornée sur $A'$. Mais $\bar c\ind_BD$ implique $\cb(\bar c/B')=\cb(\bar c/B)\in\bdd(B)$. Alors
$$\cb(\bar c/B)\in\bdd(B)\cap\bdd(A')=\bdd(A),$$
puisque $B\ind_AA'$. Donc $T$ est relativement monobasée.\ebew

\bem\label{mb-models} De même que pour le cas classique, on peut toujours se ramener au cas où $A$ et $B$ sont des modèles pour la définition de relativement monobasé sous l'hypothèse (\dag).\ebem
\bew Si $A$, $B$ et $\bar c$ témoignent que la théorie $T$ n'est pas monobasée au dessus de $T_0$ et $D\supseteq A$ est un modèle indépendant de $B\bar c$ sur $A$, alors la preuve de la proposition \ref{mb-para} montre que $D$, $\acl(BD)$ et $\bar c$ le témoignent aussi. Quant à $B$, considérons $\M\supseteq B$ un modèle indépendant de $\bar c$ sur $B$. Comme $\sscl{A\bar c}\subseteq\acl(B\bar c)$ le Lemme \ref{L:indep} implique $\M\indi i_B\sscl{A\bar c}$ pour tout $i<n$, d'où $\M\indi i_A \sscl{A\bar c}$ par transitivité. Alors $\cb(\bar c/\M)=\cb(\bar c/B)\notin\bdd(A)$.
\ebew

\bsp\label{exemples} \mbox{ }
\begin{itemize}
\item La théorie $\omega$-stable $DCF_0$  des corps différentiellement clos de caractéristique $0$ est monobasée sur la théorie $ACF_0$ des corps algébriquement clos pour la clôture algébrique (modèle-théorique) $\acl_\delta$, qui satisfait (\dag) et (\ddag) \cite[Lemme 6.12 et p.\ 111]{Po85}.
\item L'expansion d'une théorie stable $T$ sans la propriété de recouvrement fini par un prédicat unaire générique \cite{CP98} est simple et monobasée sur $T$ pour la clôture $\acl=\acl_T$ qui satisfait (\dag) et (\ddag). De même pour la modèle-complétion, si elle existe, de l'expansion de $T$ par un automorphisme \cite{CP98}, si l'on prend la clôture $\acl_\sigma$. En particulier, la théorie $ACFA$ d'un corps de différence existentiellement clos \cite{ChHr99} est supersimple, et monobasée sur la théorie $ACF$ d'un corps algébriquement clos (de même caractéristique).
\end{itemize}
\ebsp

L'hypothèse que la théorie $T$ est relativement monobasée sur ses réduits permet de montrer que le noyau de l'homomorphisme du théorème \ref{groupe-homo} est fini.

\satz\label{groupe-monobase} Soit $T$ une théorie simple avec des réduits $(T_i:i<n)$ stables qui ont l'élimination géométrique des imaginaires. Si $T$ est relativement monobasée au dessus de ses réduits pour un opérateur de clôture $\sscl.$ satisfaisant $(\dag)$ et $(\ddag)$, alors tout groupe $G$ type-définissable dans $T$ a un sous-groupe type-définissable d'indice borné qui se plonge définissablement modulo un noyau fini dans un produit fini de groupes $H_i$ interprétables dans les $T_i$.\esatz
\bew Soient $G_i$ les sous-groupes type-définissables de $G$ d'indice borné, $H_i$ les groupes $*$-interprétables dans $T_i$ et $\phi_i:G_i\to H_i$ les homomorphismes donnés par le théorème \ref{groupe-homo} pour tout $i<n$, et soit $A$ un ensemble de paramètres tel que tout soit défini sur $A$. 
Posons $N=\bigcap_{i<n}\ker(\phi_i)$, $H=\prod_{i<n}H_i$ et $\phi=\prod_{i<n} \phi_i:\bigcap_{i<n}G_i\to H$ avec $\ker(\phi)=N$.

Le théorème \ref{groupe-homo} et la remarque \ref{R:DT_i} assurent de plus que 
si $g,g'$ sont deux génériques indépendants de $\bigcap_{i<n}G_i$, il existe $D\supseteq A$ indépendant de $g,g'$ sur $A$ tel que pour tout $i<n$ on ait
$$\acl(g,D),\acl(g',D)\indi i_{\alpha_i}\acl(gg',D),$$ où
$$\alpha_i=\acl(gg'N_i,D)\cap\acl_i(\acl(gN_i,D),\acl(g'N_i,D)).$$
Posons $\alpha=(\alpha_i:i<n)$. Puisque $g\ind_Dg'$ il suit de (\dag) que
$$\sscl{\acl(g,D),\acl(g',D)}\subseteq\bigcap_{i<n}\acl_i(\acl(g,D),\acl(g',D)).$$
Or, $\alpha$ est contenu dans $\bigcup_{i<n}\acl_i(\acl(g,D),\acl(g',D))$. La condition (\ddag) entraîne l'inclusion
$$\sscl{\acl(g,D),\acl(g',D),\acl(\alpha)}\subseteq\bigcap_{i<n}\acl_i(\acl(g,D),\acl(g',D),\acl(\alpha)).$$
Ainsi, comme $\acl(\alpha)\subseteq\acl(gg',D)$, on a que pour tout $i<n$
$$\sscl{\acl(g,D),\acl(g',D),\acl(\alpha)}\indi i_{\acl(\alpha)}\acl(gg',D).$$
Cette indépendance donne 
$$\acl(g,D),\acl(g',D)\ind_{\acl(\alpha)}\acl(gg',D),$$
puisque $T$ est relativement monobasée.
On obtient que 
$$\begin{aligned}gg'&\in\acl(\alpha)=\acl(\alpha_i:i<n) \subset \acl(D,gg'N_i:i<n)\\
&=\acl(D,\phi_i(gg'):i<n)=\acl(D,\phi(gg')).\end{aligned}$$
On a ainsi vérifié que $gg'$ et $\phi(gg')$ sont interalgébriques sur $D$.

Par stabilité chaque $H_i$ est une limite projective de groupes $T_i$-interprétables et donc $H$ est une limite projective de produits de groupes $T_i$-interprétables. Comme $gg'$ est un uplet fini, il existe une projection $\pi$ de $H$ sur un produit de groupes $T_i$-interprétables $\tilde H=\prod_{i<n}\tilde H_i$ telle que $gg'$ est algébrique sur $(\pi\circ\phi)(gg'),D$. Le morphisme obtenu $\tilde\phi=\pi\circ\phi:\bigcap_{i<n}G_i\to\tilde H$ est à noyau fini.\ebew
\bem\label{super-mono} Si la composante connexe de $G$ sur $D$ est une intersection de groupes relativement définissables, alors par compacité $\tilde\phi$ est défini sur un sous-groupe de $G$ d'indice fini. Ceci est le cas si $T$ est stable ou supersimple.\ebem

Le Théorème \ref{groupe-monobase} et la Remarque \ref{super-mono} appliqués à $DCF_0$ et $ACFA$ (Exemple \ref{exemples}) permettent de retrouver le théorème de Kowalski-Pillay \cite{KP02}  et partiellement celui de Pillay \cite{Pi97} ~:
\kor\label{KoPi-Pi} Tout groupe définissable dans $DCF_0$ ou $ACFA$ a un sous-groupe d'indice fini qui se plonge définissablement dans un groupe algébrique, modulo un noyau fini.\ekor

\section{CM-trivialit\'e relative}\label{S:def}

Les théories monobasées rentrent dans un cadre plus large, les théories CM-triviales, dont on donne une définition relative qui nous permet d'étudier les groupes définissables à centre près. 

\begin{definition}\label{defcmtrivial}
La théorie $T$ est \emph{CM-triviale au dessus de $T_0$ pour $\sscl.$} si pour tous ensembles réels $A\subseteq B$ algébriquement clos et tout uplet réel $\bar c$, si 
\[ \sscl{A\bar c} \indi0_{A} B,\]
alors la base canonique $\cb(\bar c/A)$ est bornée sur $\cb(\bar c/B)$. (Notons que ces bases canoniques sont des hyperimaginaires de $T$).\newline
Plus généralement, $T$ est CM\emph{-triviale au dessus de $(T_i:i<n)$ pour $\sscl.$} si pour tout $A\subseteq B$ algébriquement clos et tout uplet $\bar c$, si
$$\sscl{A\bar c}\indi i_A B\mbox{ pour tout }i<n,$$
alors la base canonique $\cb(\bar c/A)$ est bornée sur $\cb(\bar c/B)$.\end{definition}

\begin{remark}
Toute th\'eorie est CM-triviale au dessus d'elle m\^eme pour l'opé\-ra\-teur $\acl$. Pour le même opérateur de clôture, si une th\'eorie $T$ est CM-triviale au-dessus de son réduit au langage de l'égalité, alors elle est CM-triviale au sens classique ; la réciproque est vraie si $T$ élimine géométriquement les imaginaires.

Toute théorie relativement monobasée est relativement CM-triviale. Enfin, si $T$ est stable, les bases canoniques sont imaginaires, et la clôture bornée est la clôture algébrique.\end{remark}

\begin{definition}
La th\'eorie $T$ est \emph{$2$-ample au dessus de $T_0$ pour $\sscl.$} s'il existe des uplets réels $\bar a$, $\bar b$ et $\bar c$ tels que:
\begin{itemize}
\item $ \acl(\bar a,\bar b) \indi0_{\acl(\bar a)} \sscl{\acl(\bar a),\bar c}$.
\item $ \bar c \ind_{\bar b}\bar a\bar b$.
\item $ \bar c\nind_{\bdd(\bar a)\cap\bdd(\bar b)}\bar a$.
\end{itemize}
Plus généralement, $T$ est \emph{$2$-ample au dessus de $(T_i:i<n)$ pour $\sscl.$} s'il existe des uplets $\bar a$, $\bar b$ et $\bar c$ tels que:
\begin{itemize}
\item $ \acl(\bar a,\bar b) \indi i_{\acl(\bar a)} \sscl{\acl(\bar a),\bar c}$ pour tout $i<n$.
\item $ \bar c \ind_{\bar b}\bar a\bar b$.
\item $ \bar c\nind_{\bdd(\bar a)\cap\bdd(\bar b)}\bar a$.
\end{itemize}
\end{definition}
\satzli\label{L:CM=2ample}
Une th\'eorie $T$ est CM-triviale au-dessus de $(T_i:i<n)$ si et seulement si elle n'est pas $2$-ample au dessus de $(T_i:i<n)$.\esatzli
\bew
Supposons que $T$ soit $2$-ample pour $\sscl.$, propriété témoignée par des uplets $\bar a$, $\bar b$ et $\bar c$. Posons  $A=\acl(\bar a)$ et $B=\acl(\bar a,\bar b)$.
Par d\'efinition $\sscl{A\bar c} \indi i_A B$ pour tout $i<n$.
Puisque $\bar c\nind_{\bdd(\bar a)\cap\bdd(\bar b)}\bar a$, la base canonique $\cb(\bar c/A)$ n'est pas bornée sur $\bar b$. D'autre part, $\bar c\ind_{\bar b} B$ implique $\cb(\bar c/B)\in\bdd(\bar b)$, d'o\`u $\cb(\bar c/A)$ n'est pas bornée sur $\cb(\bar c/B)$. Donc $T$ n'est pas CM-triviale au-dessus de $(T_i:i<n)$.

Réciproquement, supposons que $T$ ne soit pas CM-triviale au-dessus de $(T_i:i<n)$ pour $\sscl.$, propriété t\'emoignée par des ensembles $A$, $B$ et un uplet $\bar c$. On pose $\bar a=\sscl{\bar a}=A$ et $\beta= \cb(\bar c/B)$. Prenons un modèle $M$ tel que $M\ind_\beta B\bar c$ et $\beta\in\bdd(M)$. 
Considérons un uplet réel $\bar b$ dans $M$ qui algébrise $\beta$. Par définition de la base canonique, $\bar c\ind_\beta B$ donc $\bar c\ind_\beta M B$.  Alors $\bar c\ind_{\bar b} \bar a \bar b$. 
De plus, $\beta\in\bdd(B)$ implique $\bar c \ind_B M$, ce qui entraîne $\acl(\bar a\bar b)\ind_B \sscl{\bar a\bar c}$ car $\bar a \in B$ et $\sscl{\bar a\bar c} \subset \acl(\bar a \bar c)$. Alors $\acl(\bar a\bar b)\indi i_B \sscl{\bar a\bar c}$. Par hypothèse $B\indi i_{\bar a} \sscl{\bar a \bar c}$~; donc par transitivité, on conclut   $\acl(\bar a\bar b)\indi i_A\sscl{A\bar c}$ pour tout $i<n$. Enfin,
$\cb(\bar c/A)\notin\bdd(\bar b)$, car $\bar b \ind_\beta \bar a$ et $\cb(\bar c/A)\notin\bdd(\beta)$. Donc, $\bar c\nind_{\bdd(\bar a)\cap \bdd(\bar b)}\bar a$.\ebew

\satzli Supposons que $\sscl.$ satisfasse $(\dag)$. Alors la CM-trivialité relative est conservée par l'adjonction ou la suppression de paramètres. \esatzli
\bew La conservation par adjonction de paramètres est évidente. Pour la suppression, supposons que $D$ soit un ensemble de paramètres tel que $T(D)$ soit CM-triviale au dessus de $(T_i:i<n)$, et considérons $A\subseteq B$ algébriquement clos, ainsi qu'un uplet $\bar c$, avec $\sscl{A\bar c}\indi i_AB$ pour tout $i<n$. On peut toujours supposer $AB\bar c\ind D$. Soit $A'=\acl(AD)$ et $B'=\acl(BD)$. Alors $B\bar c\ind_A D$, et comme dans la preuve de la proposition \ref{mb-para} on obtient $\sscl{A'\bar c}\indi i_{A'}B'$ pour tout $i<n$.

Par hypothèse, $\cb(\bar c/A')\in\bdd(\cb(\bar c/B'),D)$. Or, $\bar c\ind_AD$, d'o\`u $\cb(\bar c/A')=\cb(\bar c/A)\in\bdd(A)$~; de même $\bar c\ind_BD$ donne $\cb(\bar c/B')=\cb(\bar c/B)\in\bdd(B)$.
Alors $AB\ind D$ implique $\cb(\bar c/A)\ind_{\cb(\bar c/B)}D$, et donc $\cb(\bar c/A)\in\bdd(\cb(\bar c/B))$.
\ebew

\begin{remark}\label{models} Comme dans le cas classique \cite[Corollary 2.5]{Pi95}, sous l'hypothèse (\dag) on peut toujours se ramener au cas où $A$ et $B$ sont des modèles dans la définition de CM-trivialité relative.\end{remark}
\bew Soit $D\supseteq A$ un modèle indépendant de $B\bar c$ sur $A$, et $B'=\acl(BD)$. Alors comme dans la preuve de la proposition \ref{mb-para} on a $\sscl{D\bar c}\indi i_D B'$ pour tout $i<n$. De plus, $\cb(\bar c/D)=\cb(\bar c/A)$ et $\cb(\bar c/B')=\cb(\bar c/B)$. On peut donc remplacer $A$ par $D$.

Quant à $B$, considérons $\M\supseteq B$ un modèle indépendant de $\bar c$ sur $B$. Alors $\sscl{A\bar c}\indi i_A\M$ pour $i<n$ comme dans la preuve de la remarque \ref{mb-models}, et $\cb(\bar c/\M)=\cb(\bar c/B)$. On peut donc remplacer $B$ par $\M$.\ebew

Rappelons que les groupes définissables dans une structure CM-triviale de rang de Morley fini sont  nilpotent-par-finis car une telle structure n'interprète ni de corps infini, ni de mauvais groupes \cite{Pi95}. Pour le cas relativement CM-trivial, nous allons montrer qu'un groupe définissable peut être vu  modulo un noyau central comme un sous-groupe d'un produit de groupes interprétables dans les réduits. Si $T$ est stable, le centre d'un groupe est relativement définissable par la condition de chaîne sur les centralisateurs. Par contre, pour un groupe $G$ type-définissable dans une théorie simple, son centre ne l'est pas nécessairement. On utilisera alors le centre approximatif \cite[Definition 4.4.9]{wa00}
$$\tilde Z(G)=\{g\in G:[G:C_G(g)]<\infty\}.$$
On note que d'après \cite[Lemma 4.2.6]{wa00} le centre approximatif est relativement définissable, et que pour $G$ stable et connexe, $\tilde Z(G) = Z(G)$. De plus, d'après \cite[Proposition 4.4.10]{wa00}, si $G^0$ est la composante connexe de $G$ sur les paramètres qui servent à définir $G$, alors $[\tilde Z(G),G^0]$ est un sous-groupe fini central de $G^0$, modulo lequel $\tilde Z(G^0)=\tilde Z(G)\cap G^0$ est central dans $G^0$.

\satz\label{T:homo} Soit $T$ une théorie simple avec des réduits stables $(T_i:i<n)$ qui ont l'élimination géométrique des imaginaires. Si $T$ est relativement CM-triviale au-dessus de ces réduits pour un opérateur de clôture $\sscl.$ satisfaisant $(\dag)$ et $(\ddag)$, alors tout groupe $G$ type-définissable dans $T$ a un sous-groupe type-définissable d'indice borné qui s'envoie définissablement dans un produit fini de groupes interprétables dans les $T_i$, où le noyau du morphisme est contenu (à indice fini près) dans le centre approximatif $\tilde Z(G)$ de $G$.\esatz

\bew Soient $G_i$ les sous-groupes type-définissables de $G$ d'indice borné, $H_i$ les groupes $*$-interprétables dans $T_i$ et $\phi_i:G_i\to H_i$ les homomorphismes donnés par le théorème \ref{groupe-homo} pour tout $i<n$. On ajoute les paramètres nécessaires au langage, et on remplace par la suite $G$ par $\bigcap_{i<n}G_i$. 
Notons $N_i$ le noyau de $\phi_i$ pour chaque $i$ et posons $N=\bigcap_{i<n}N_i$, $H=\prod_{i<n}H_i$ et $\phi=\prod_{i<n} \phi_i:G\to H$ avec $\ker(\phi)=N$.

Le théorème \ref{groupe-homo} et la remarque \ref{R:DT_i} entraînent de plus que pour $g,g'$ deux génériques indépendants de $G$ il existe un ensemble $D$ d'éléments génériques indépendants sur $g,g'$
tel que pour $$\alpha_i=\acl(gg'N_i,D)\cap\acl_i(\acl(gN_i,D),\acl(g'N_i,D))$$
on ait
$$\acl(g,D),\acl(g',D)\indi i_{\alpha_i}\acl(gg',D)$$
pour tout $i<n$. 

Posons $\alpha = (\alpha_i:i<n)$. Comme $\alpha$ est algébrique sur $D,gg'$, on a
$$\acl(D,g),\acl(D,g')\indi i_{\acl(D,\alpha)}\acl(D,gg').$$
Puisque $\sscl.$ satisfait $(\dag)$ et $g \ind_D g'$, on a 
$$\sscl{\acl(D,g),\acl(D,g')}\indi i_{\acl(D,\alpha)}\acl(D,gg').$$
Or, $\alpha\subseteq\bigcup_{j<n}\acl_j(\acl(D,g),\acl(D,g'))$, et  par $(\ddag)$
$$\sscl{\acl(D,g),\acl(D,g'),\acl(D,\alpha)} \subseteq \acl_i(\sscl{\acl(D,g),\acl(D,g')},\acl(D,\alpha)) $$ donc
$$\sscl{\acl(D,g),\acl(D,g'),\acl(D,\alpha)}\indi i_{\acl(D,\alpha)}\acl(D,gg').$$
Comme $\phi(gg')$ et $\alpha$ sont interalgébriques sur $D$ d'après la remarque \ref{homo-prec} et la définissabilité de $\phi$, on obtient l'indépendance 
\begin{equation}\label{eqn2}\sscl{\acl(D,g),\acl(D,g'),\acl(D,\phi(gg'))}\indi i_{\acl(D,\phi(gg'))}\acl(D,gg').\end{equation}

Soient $a,b,e\in D$ distincts, et $D'=D\setminus\{e\}$. On pose
$$A=\acl(D',\phi(gg')),\quad B=\acl(D',gg')\quad\mbox{et}\quad\bar c=(g,g',age,e^{-1}g'b).$$
Comme $age\in\acl(D,g)$ et $e^{-1}g'b\in\acl(D,g')$, on obtient
$$\sscl{g,g',age, e^{-1}g'b,\acl(D',\phi(gg'))}\indi i_{\acl(D',\phi(gg'),e)}\acl(D',gg').$$
Puisque $e\ind_{\acl(D',\phi(gg'))}gg'$, on a par le Lemme \ref{L:indep}
$$\acl(D',\phi(gg'),e)\indi i_{\acl(D',\phi(gg'))}\acl(D',gg')$$
et donc par transitivité
$$\sscl{g,g',age,e^{-1}g'b,\acl(D',\phi(gg'))}\indi i_{\acl(D',\phi(gg'))}\acl(D',gg').$$
Autrement dit, pour tout $i<n$
$$\sscl{\bar c,A} \indi i_{A} B.$$
Par CM-trivialité relative, $\cb(\bar c/A)\in\bdd(\cb(\bar c/B)$.

Soit $Z$ le centralisateur approximatif de $N$ dans $G$, c'est-à-dire
$$Z=\tilde C_G(N)=\{g\in G:[N:C_N(g)]<\infty\}.$$
(Si $T$ est stable, le centralisateur $C_G(N)$ est relativement définissable et il peut remplacer $Z$ dans la suite.) D'après \cite[Lemma 4.2.6]{wa00} le centralisateur approximatif est un sous-groupe relativement définissable, disons $Z=G\cap X$ pour un ensemble $\emptyset$-définissable $X$. Alors 
$$(x^{-1}y\in X\land x\in G\land y\in G)\lor x=y$$
est une relation d'équivalence type-définissable, et le translaté $aZ$ en est la classe de $a$~: il s'agit bien d'un élément imaginaire qui jouera un rôle important dans le lemme suivant qui nous permettra de conclure.
\lmm Si $N\cap \tilde Z(G)$ n'est pas d'indice fini dans $N$, alors $aZ\in\acl^{eq}(\cb(\bar c/A))\setminus\acl^{eq}(\cb(\bar c/B))$.\elmm
\bew Notons d'abord que puisque $\tilde Z(G)$ est relativement définissable, l'indice de $N\cap\tilde Z(G)$ dans $N$ est fini si et seulement s'il est borné. Supposons que $N\cap \tilde Z(G)$ n'est pas d'indice fini dans $N$. Alors $Z$ n'est pas d'indice fini dans $G$~: sinon, il contiendrait un générique $g$ de $G$. Par définition, $N$ aurait un générique $n$ indépendant de $g$ avec $[g,n]=1$. Cela entraînerait que $[G:C_G(n)]<\infty$ et donc $n$ serait dans $\tilde Z(G) \cap N$, ce qui impliquerait que l'indice $\tilde Z(G) \cap N$ dans $N$ fût fini. En particulier le translaté $aZ$ n'est pas dans $\bdd(\emptyset)$. 

Calculons la base canonique $\cb(\bar c/B)$. Puisque $g$ et $age$ sont deux génériques indépendants sur $(D',gg')$, l'élément $agg'b$ est définissable sur $(D',gg')$ et $\bar c$ est interalgébrique sur $(gg',agg'b)$ avec $(g,age)$, on obtient que $$\bar c\ind_{gg',agg'b}D'$$ et $B\cap\acl(\bar c)=\acl(gg',agg'b)$.
En particulier $\cb(\bar c/B)$ est interalgébrique avec $gg',agg'b$. Or, l'indépendance $$a\ind gg',agg'b$$ implique que $aZ\ind\cb(\bar c/B)$, et donc $aZ$ n'est pas algébrique sur $\cb(\bar c/B)$.

Montrons maintenant que $aZ$ est algébrique sur $\cb(\bar c/A)$. Par définition
$$\bar c \ind_{\cb(\bar c/A)} A$$
et donc pour une suite de Morley $J$ de $\lstp(\bar c/A)$, on a 
$$J \ind_{\cb(\bar c/A)} A.$$ Comme $aZ$ est algébrique sur $a \in A$ il suffit de voir que $aZ$ l'est sur $J$. Pour la suite notons $$J=(g_i,g'_i,ag_ie_i, e_i^{-1}g'_ib:i\in I)$$ et 
considérons une autre réalisation $a'$ de $\lstp(a/J)$. Soit $b'$ tel que $(a',b')\models\lstp(a,b/J)$. Alors $ag_ig'_ib=a'g_ig'_ib'$ pour tout $i\in I$ car ce produit est définissable sur $J$.

Comme $gg'\ind_{\phi(gg')}A$, la suite $(g_ig'_i:i\in I)$ est indépendante sur $\phi(gg')$. Si la suite est suffisamment longue, il existe un élément $g_ig'_i$ indépendant de $a,b,a',b',gg'$ sur $\phi(gg')$~; cet élément est donc générique dans  $gg'N$ sur $a,b,a',b',gg'$. Alors $g_ig'_i = gg'n$ pour $n$ générique dans $N$ sur  $gg',a,b,a',b'$. Il suit que
$$n^{-1}(gg')^{-1}a^{\prime-1}agg'n=b'b^{-1}.$$
Prenons $m \models \lstp(n/gg',a,b,a',b')$ indépendant de $n$. Alors
$$m^{-1}(gg')^{-1}a^{\prime-1}agg'm=b'b^{-1}=n^{-1}(gg')^{-1}a^{\prime-1}agg'n$$
et 
$$(mn^{-1})^{-1}(gg')^{-1}a^{\prime-1}agg'(mn^{-1})=(gg')^{-1}a^{\prime-1}agg'.$$
Comme  $mn^{-1}$ est générique dans $N$ sur $gg',a,b,a',b'$, l'élément  $(gg')^{-1}a^{\prime-1}agg'$ commute avec un générique et donc appartient à $\tilde C_G(N)=Z$. La normalité de  $N$ dans $G$ entraîne celle de $Z$.  Donc $a^{\prime-1}a$ est dans $Z$, d'où  $aZ=a'Z$ est algébrique sur $J$.\ebew

Comme $\cb(\bar c/A)\in\bdd(\cb(\bar c/B))$, le sous-groupe  $N\cap\tilde Z(G)$ est d'indice fini dans $N$; en particulier l'élément $gg'\tilde Z(G)$ est algébrique sur $\phi(gg')$. Du fait que $\tilde Z(G)$ est relativement définissable dans $G$, l'imaginaire fini $gg'\tilde Z(G)$ est donc algébrique sur une partie finie de $\phi(gg')$. La stabilité des réduits permet comme à la fin de la preuve du théorème \ref{groupe-monobase} de trouver une projection $\pi$ de $H$ sur un produit $\tilde H=\prod_{i<n}\tilde H_i$ de groupes $T_i$-interprétables telle que  $gg'\tilde Z(G)$ reste algébrique sur $(\pi\circ\phi)(gg')$. Le morphisme $\tilde\phi=\pi\circ\phi : G\to\tilde H$ a alors un noyau contenu (à indice fini près) dans $\tilde Z(G)$.\ebew

\bem Si la composante connexe de $G$ est une intersection de groupes relativement définissables, alors par compacité $\tilde\phi$ est défini sur un sous-groupe de $G$ d'indice fini.\ebem

\kor\label{cm-simple} Avec les hypothèses du théorème précédent, tout groupe $G$ simple type-définissable dans $T$ se plonge définissablement dans un groupe interprétable dans l'un des réduits.\ekor
\bew Par le théorème \ref{T:homo} il y a un homomorphisme non-trivial $\phi=\prod_{i<n}\phi_i$ de $G$ dans un produit $H=\prod_{i<n}H_i$ de groupes interprétables dans les réduits. Puisque $G$ est simple, il se plonge dans l'un des $H_i$. \ebew

\kor\label{C:corps} Avec les hypothèses du théorème précédent, tout corps $K$ type-définissable dans $T$ est définissablement isomorphe à un sous-corps d'un corps interprétable dans l'un des réduits. Si $T$ est de rang $SU$ fini alors $K$ est algébriquement clos et définissablement isomorphe à un corps interprétable dans un des réduits.\ekor

\bew Considérons le groupe algébrique simple $PSL_2(K)$. Par le corollaire précédent, il se plonge définissablement dans un groupe $H$ interprétable dans l'un des réduits $T_i$. Ce plongement nous permet de supposer que $K^+\rtimes K^\times$ est un sous-groupe de $H$, sachant que l'action par translation de $K^\times$ sur $ K^+$ se traduit dans $H$ par l'action par conjugaison. Soient $\bar K^+$ et $\bar K^\times$ les plus petits sous-groupes $T_i$-type-interprétables de $H$ contenant $K^+$ et $K^\times$ respectivement. 
Notons que $K^\times$ normalise aussi $\bar K^+$, d'où $\bar K^\times\le N_H(\bar K^+)$. Donc $\bar K^\times$ agit par conjugaison sur $\bar K^+$. Or, d'après \cite[Remark 27]{fW90} ou \cite[Theorem 2.1.19]{Wa97} les groupes $\bar K^+$ et $\bar K^\times$ satisfont les mêmes propriétés génériques que $K^+$ et $K^\times$. 
Donc pour tout  $a\in\bar K^+$ générique il y a $b\in\bar K^\times$ avec $b\cdot 1_{K^+}=a$. Comme $\bar K^\times$ est abélien, tout élément générique $a\in\bar K^+$ a le même fixateur $F=C_{\bar K^\times}(a)$. Tout élément de $\bar K^+$ étant la somme de deux génériques, $F$ est aussi le fixateur de $\bar K^+$ tout entier, $F=C_{\bar K^\times}(\bar K^+)$. Comme $F\cap K^\times=\{1\}$ on peut quotienter par $F$ et supposer que $\bar K^\times$ ne fixe aucun élément générique de $\bar K^+$.

On peut identifier génériquement $\bar K^+$ et $\bar K^\times$ en associant à $a\in\bar K^+$ générique l'unique $b\in\bar K^\times$ avec $a = b\cdot 1_{K^+}$. Ainsi on définit génériquement un produit distributif sur $\bar K^+$, qu'on étend à $\bar K^+$ tout entier~: pour $a,b\in\bar K^+$ on trouve $a_1,a_2,b_1,b_2\in\bar K^+$ génériques tels que $a=a_1+a_2$ et $b=b_1+b_2$, et chaque $a_i$ est indépendant de chaque $b_j$. On pose alors
$$a\cdot b=(a_1+a_2)(b_1+b_2)=a_1b_1+a_1b_2+a_2b_1+a_2b_2~.$$
On vérifie que la définition ne dépend pas du choix des décompositions : en effet si $b$ est générique et si $a_1$, $a_2$, $a'_1$ et $a'_2$ sont génériques sur $b$ tels que $a_1+a_2 = a'_1 +a'_2$ et $a_2$ est générique sur $a'_1$, $a'_2$ et $b$ alors la distributivité générique entraîne 
$$a_1b = (a'_1 + (a'_2-a_2))b = a'_1 b +  (a'_2-a_2) b = a'_1 b +a'_2 b-a_2 b.$$ 

Ainsi on obtient une structure de corps $T_i$-type-interprétable sur $\bar K^+$ qui étend celle de $K$. Or, tout corps type-interprétable dans une théorie stable est contenu dans un corps interprétable d'après \cite[Corollaire 5.20]{Po87}.

Le plongement de $K$ dans $\bar K$ nous donne une paire de corps infinis interprétables dans $T$. Si $T$ est de rang $SU$ fini, l'extension $\bar K$ sur $K$ est nécessairement de degré fini.
Or, $T_i$ est superstable de rang borné par le rang $SU$ de $T$. Donc $\bar K$ est algébriquement clos \cite{Mac71,ChSh80}. Comme le degré de l'extension est fini, le corps $K$ est réel clos ou algébriquement clos. La simplicité de $T$ interdit le premier cas et donc les deux corps sont égaux.\ebew

\bem\label{cm-homo} Sans la condition $(\ddag)$ on peut montrer que l'homomorphisme du théo\-rème \ref{T:homo} n'est pas trivial si $G$ est non-isogène à un groupe abélien~: si $N$ était d'indice borné dans $G$, la classe $gg'N$ serait dans $\bdd(\emptyset)$, et $\alpha_i$ et $\phi(gg')$ seraient algébriques sur $\emptyset$. L'indépendance 
$$\acl(D,g),\acl(D,g')\indi i_{\alpha_i}\acl(D,gg')$$
donnerait trivialement, en utilisant (\dag), l'équation (\ref{eqn2}) de la démonstration du theorème \ref{T:homo}; la suite de cette preuve n'utilise plus la condition $(\ddag)$ et montre que $N$ est contenu, à indice fini près, dans $\tilde Z(G)$.

Par conséquent les deux corollaires du théorème restent vrais sans l'hypothèse~$(\ddag)$.\ebem

\section{Amalgames}\label{S:Amal}
Dans cette partie on montre la CM-trivialité relative des corps colorés dans \cite{Po99,Po01,BaHo,BMPZ04,BMPZ05,BHPW06} et des fusions de deux théories  dans \cite{Hr92,BMPZ06,Z08}. Ces structures sont obtenues par amalgamation à la Hrushovski-Fraïssé (libre ou collapsée) à l'aide d'une certaine fonction, la {\em prédimension}, qui induit un nouvel opérateur, \emph{la clôture autosuffisante}. Notre approche consiste à isoler des caractéristiques communes de ces exemples pour montrer de façon uniforme la CM-trivialité relative au-dessus des théories de base par rapport à la clôture autosuffisante. Même si Hrushovski a déjà énoncé une caractérisation des groupes définissables pour la fusion sur l'égalité, le théorème \ref{T:homo} et ses corollaires nous permettent en particulier d'obtenir de nouveaux résultats sur la définissabilité de corps et de mauvais groupes pour les corps colorés et les fusions.

\bem Les constructions {\em ab initio}
dans un langage relationnel peuvent être vues comme des expansions colorées de la théorie de l'égalité, où des $n$-uplets sont colorés s'ils sont liés par une des relations. Néanmoins dans ce cas la CM-trivialité relative n'est qu'un petit renforcement de la CM-trivialité absolue qui est déjà connue \cite{Yo03,VY03}. De toute façon, il n'y a pas de groupes définissables infinis.

Le \emph{groupe de Baudisch} \cite{Ba96} est un
$\mathbb{F}_p$-espace vectoriel $V$ muni d'une sorte additionnelle pour le produit extérieur ${\bigwedge}^2V$ avec la forme bilinéaire naturelle.
Baudisch a montré que ce groupe est CM-trivial et donc que tout groupe définissable est nilpotent-par-fini.
Nous conjecturons que ce groupe (collapsé ou non) est relativement CM-trivial au-dessus de l'espace vectoriel $V$.
Ceci permettrait de conclure que la classe de nilpotence de tout groupe $G$ connexe définissable soit au plus $2$.

En revanche le corps différentiel rouge \cite{BMP09} sort de notre  cadre puisque la clôture autosuffisante n'est pas contenue dans la clôture algébrique.\ebem

Nous revenons maintenant à une description succincte des constructions par amalgamation des corps colorés et des fusions, afin de mettre en évidence les caractéris\-tiques qui nous sont utiles pour vérifier la propriété de CM-trivialité relative~:
\begin{itemize}
\item[\bf Coloré] On se donne une théorie de base $T_0$ et un nouveau prédicat $P$, dont les points sont dits {\em colorés}, et on considère la classe $\mathcal F$ des modèles colorés de $T_0^\forall$ ; 
\item[\bf Fusion] On considère plusieurs théories $T_i$  avec un réduit commun $T_{com}$, et $\mathcal F$ dénote la classe des modèles de chaque $T_i^\forall$. 
\end{itemize}

La méthode d'amalgamation consiste à construire une structure avec une géométrie prédéterminée. Pour cela, on introduit une prédimension $\delta$  sur les modèles finiment engendrés de $\mathcal F$ qui satisfait l'inégalité sous-modulaire 
$$\delta(A\cup B)\le\delta(A)+\delta(B)-\delta(A\cap B).$$
(Par abus de notation on écrit $\delta(A)$ pour la prédimension de la structure engendrée par $A$.)

La fusion de deux théories $T_1$ et $T_2$ de rang de Morley fini et de rang et degré définissables sur l'égalité se sert de la prédimension 
$$\delta(A) = n_1\RM_1(A) + n_2\RM_2(A)-n|A|,$$
où $n=n_1\RM(T_1)=n_2\RM(T_2)$.

La fusion de deux théories fortement minimales sur un  $\F_p$-espace vectoriel infini utilise la prédimension
$$\delta(A) = \RM_1(A) + \RM_2(A)-\ldim_{\F_p}(A).$$

Les corps colorés sont des corps algébriquement clos avec un prédicat $P$ pour un sous-ensemble, avec prédimension
$$\delta(k)=2 \, \mathrm{degtr}(k) - \dim_P(P(k)).$$
Dans le corps {\em noir}, le prédicat $P$ dénote un sous-ensemble $N$ et $\dim_P(N)=|N|$. Le corps {\em rouge} est de caractéristique positive $p$ et $P$ dénote un sous-groupe additif propre $R$ avec $\dim_P(R)=\ldim_{\F_p}(R)$. Enfin, le corps {\em vert} est de caractéristique $0$ et $P$ dénote un sous-groupe multiplicatif divisible sans torsion $\ddot{U}$ (vu comme un $\Q$-espace vectoriel) avec $\dim_P(\ddot U)=\ldim_{\Q}(\ddot{U})$.

En travaillant sur une sous-structure $A$, on peut aussi définir une version relative $\delta({.}/A)$ de la prédimension. La sous-modularité implique 
$$\delta(\bar a/A)=\lim_{A_0\to A}\{\delta(\bar a/A_0)\}=\inf\{\delta(\bar a\cup A_0)-\delta(A_0): A_0 \subseteq A \text{ finiment engendrée}\}.$$
En particulier, la limite existe et la prédimension relative est également sous-modulaire. (Par le même abus de notation
nous écrirons $\delta({.}/A)$ quand il s'agit de la prédimension au-dessus de la structure engendrée par $A$.)

Selon la partie négative de la prédimension on distingue les deux comportements suivants :
\begin{itemize}
\item[{\bf Dégénéré}] La partie négative de la prédimension correspond à la cardinalité d'un certain prédicat (ou, plus généralement, à la dimension d'une pré\-géo\-mé\-trie dégénérée). Par exemple, la fusion de deux théories sur l'égalité ou le corps noir.
\item[{\bf Modulaire}] Il y a un groupe abélien $\emptyset$-définissable dans le langage commun à toutes les théories qui est soit un $\mathbb F_p$-espace vectoriel (le cas de la fusion de deux théories sur un $\mathbb F_p$-espace vectoriel ou des corps rouges), soit un groupe divisible avec $n$-torsion finie pour chaque $n$ (le cas des corps verts). Toute structure est munie de cette loi de groupe et les points colorés (s'il en existe) en forment un sous-groupe (divisible si le groupe l'est).\end{itemize}

Pour le cas modulaire, on supposera pour simplifier la rédaction que le groupe abélien $\emptyset$-définissable a même domaine que la structure (pour le corps vert où il s'agit du groupe multiplicatif, cette hypothèse est vérifiée si on met de coté l'élément $0$). Grâce à une Morleyisation, on suppose pour la suite que les théories de départ ($T_i$, $T_{com}$) 
éliminent les quanteurs et que le langage est relationnel sauf pour la loi de groupe dans le cas modulaire (c'est-à-dire l'addition dans le corps rouge, la multiplication dans le corps vert, et l'addition vectorielle dans la fusion sur un espace vectoriel).
Si ce groupe est divisible avec $n$-torsion finie pour chaque $n$, par abus de langage, on entendra par structure engendrée par une partie $B$, la clôture divisible du groupe engendré par $B$. Ainsi une structure finiment engendrée correspondra à un groupe divisible de rang fini (au sens groupe-théorique).
 
Par élimination des quanteurs dans $T_i$, le $i$-diagramme d'une structure $A$ détermine son $i$-type, et\begin{itemize}
\item[$(\delta)$] $\quad\bigcup_i\diag_i(A)$, plus la coloration de $A$ (s'il y en a), détermine $\delta(A)$.\end{itemize}

On se restreint à la sous-classe $\K$ des structures de $\mathcal F$ telles que toute sous-structure finiment engendrée a une prédimension positive ou nulle.

Pour $M$ dans $\K$, une sous-structure $A$ est {\em autosuffisante} dans $M$, noté $A\le M$, si $\delta(\bar a/A)\ge0$ pour tout uplet fini $\bar a\in M$. La sous-modularité entraîne que l'intersection de deux structures autosuffisantes l'est aussi, et ainsi l'existence pour tout $A\subseteq M$ d'une plus petite structure autosuffisante contenant $A$. On l'appelle la {\em clôture autosuffisante} de $A$ et on la note $\sscl{A}_M$; si $M\le N$ alors $\sscl{A}_M=\sscl{A}_N$ et nous noterons simplement $\sscl{A}$. Comme $\sscl{A}$ 
est unique, elle est contenue dans $\acl(A)$. Encore par sous-modularité
$$\sscl{A}=\bigcup\{\sscl{A_0}:A_0\subset A\mbox{ finiment engendrée}\}$$
et
$$\delta(\sscl{\bar a})=\liminf\{\delta(\bar a_0):\bar a_0\supseteq\bar a\mbox{ fini}\}.$$
On note que si $\delta$ prend uniquement des valeurs entières, la clôture autosuffisante d'un ensemble finiment engendré est lui-même finiment engendré. Ceci est notamment satisfait dans les exemples cités ci-dessus.

\bem\label{sscl} Si un uplet fini $\bar a$ est tel que $\delta(\bar a/A)<0$ mais $\delta(\bar a'/A)\ge0$ pour tout sous-uplet strict $\bar a'$ de $\bar a$, alors $\bar a\in\sscl{A}$.\ebem
\bew Supposons que $\bar a' = \sscl{A} \cap \bar a$ est une partie propre de $\bar a$. La sous-modularité donne alors 
$ \delta(\bar a/\sscl{A}) \leq  \delta(\bar a / A \cup \bar a')<0$.\ebew

\`A partir d'une sous-classe $\K_0\subseteq\K$ de structures finiment engendrées avec la propriété d'amalgamation pour les plongements autosuffisants, on construit par la méthode de Fraïssé un amalgame dénombrable $\M$ universel et fortement homogène pour les sous-structures autosuffisantes. Le modèle générique $\M$ s'avère être un modèle de $\bigcup_{i<n} T_i$; il est stable, superstable si $\delta$ ne prend que des valeurs entières et $\omega$-stable si on peut borner les multiplicités des formules \cite{Wa10}. Il est  de rang fini dans le cas collapsé (obtenu par le choix d'une sous-classe de $\K$ suffisamment restreinte pour rendre algébriques (certains) éléments de dimension $0$).

L'indépendance au sens de $T$ est caractérisée de la manière suivante~:\smallskip

$\mathbf{(\star)}\quad$ Pour deux uplets $\bar a$, $\bar b$ et un ensemble $C$ algébriquement clos, 
$$ \bar a\ind_C \bar b$$
\begin{center}si et seulement si\end{center}\vskip2mm
$$\parbox{7em}{$\displaystyle\sscl{\bar a\cup C}\indi i_C \sscl{\bar b\cup C}$\newline pour tout $i$} \quad\mbox{et}\quad\left\{\parbox{23em}{\begin{itemize}
\item[{\bf Dégénéré}] $\sscl{\bar a\bar b\cup C}=\sscl{\bar a\cup C}\cup\sscl{\bar b\cup C}. $
\item[{\bf Modulaire}] $\sscl{\bar a\bar b\cup C}$ est  le sous-groupe engendré par $\sscl{\bar a\cup C}$ et $\sscl{\bar b\cup C}$, et ses points colorés (s'il y en a) sont les produits de ceux de $\sscl{\bar a\cup C}$ et de $\sscl{\bar b\cup C}$.
\end{itemize}}\right.$$
\smallskip

Remarquons que les prédimensions considérées impliquent que \begin{itemize}
\item[$(\sstar)_1$]$\quad \delta(\bar a/A)\le 0$ pour tout $\bar a\in\bigcup_i\acl_i(A)$.\end{itemize}
Dans le cas coloré, il y a un unique réduit et la partie positive de la prédimension est alors nulle si $\bar a\in  \acl_0(A)$.
Pour les fusions, chacune des parties positives de la prédimension est majorée par la partie négative. 

De plus, pour les classes d'amalgamation $\K_0$ considérées, l'amalgame libre de deux  structures $B$ et $C$ de $\K_0$ au-dessus d'une sous-structure autosuffisante $A$ reste dans $\K_0$ du moment que la prédimension de $B$ sur $A$ est strictement positive et $A$ est autosuffisante maximale dans $B$. Ceci implique que
\begin{itemize}
\item[$(\sstar)_2$]$\quad$ pour tout uplet fini $\bar a\in\acl(A)$ il y a un uplet fini $\bar b$ tel que $\bar a\subset \bar b\in\acl(A)$ et $\delta(\bar b/\sscl A)=0$.\end{itemize}

\lmm La propriété $(\star)$ implique $(\dag)$ et  $(\sstar)_1+(\sstar)_2$ implique $(\ddag)$ pour la clôture autosuffisante.\elmm
\bew La première implication est évidente dans le cas dégénéré et elle suit dans le cas modulaire du fait que la loi du groupe est définie dans le langage commun.

Pour la seconde implication, on considère un ensemble $A$ autosuffisant et un uplet fini $\bar c\in\bigcup_i\acl_i(A)$. Il suffit alors de montrer que la structure engendrée par $A \cup B$ est autosuffisante où $B=\acl(\bar c)$.

On montre qu'elle est réunion croissante de sous-structures autosuffisantes dont chacune est clôture autosuffisante d'un uplet fini. Soit $\bar a$ un uplet fini de $A$ tel $\bar c\in\bigcup_i\acl_i(\bar a)$. Soit $\bar b$ un  uplet fini de  $B$ contenant $\bar c$. La condition $(\sstar)_2$ permet de supposer de plus que $\delta(\bar b/\sscl{\bar c})=0$ et donc que  $\delta(\sscl{\bar b}) = \delta(\sscl{\bar c})$.  De plus $(\sstar)_1$ implique que $0 \leq  \delta(\bar c/\sscl{\bar a}) \leq 0$. La structure $C$ engendrée par $\bar c \cup \sscl{\bar a}$ est donc autosuffisante. Donc $\sscl{\bar c}\subseteq C\cap \sscl{\bar b}$. 
Montrons que la structure engendrée par $\sscl{\bar a} \cup \sscl{\bar b}$ est autosuffisante. Comme cette structure est égale à celle engendrée par $C \cup \sscl{\bar b}$, il suffit de montrer que $\delta(\sscl{\bar b}/C) = 0$.
Par sous-modularité, 
$$\delta(\sscl{\bar b}/C) \le \delta(\sscl{\bar b}/ C \cap \sscl{\bar b}).$$
Comme $$\delta(C \cap \sscl{\bar b})\ge\delta(\sscl{\bar c})=\delta(\sscl{\bar b}),$$
on conclut que $\delta(\sscl{\bar b}/C) = 0$. \ebew

Nous pouvons maintenant montrer la CM-trivialité relative des amalgames considérés.
\satz\label{T:amalgames} Soit $\M$ l'un des modèles génériques construits dans \cite{Hr92,Po99,Po01,BaHo,BMPZ04,BMPZ05,BMPZ06,BHPW06,Z08} (de domaine $M$). Alors sa théorie $T$ est  CM-triviale pour l'opérateur de clôture autosuffisante $\sscl.$ au-dessus des théories de base $(T_i:i<n)$ associées.\esatz
\bem Notons  que les théories de base dans la fusion \cite{Z08} n'éliminent pas toujours géométriquement les imaginaires~; en revanche le Lemme \ref{L:CM=2ample} n'utilise pas cette hypothèse et  la conclusion  du Lemme \ref{L:indep} suit de la caractérisation de l'indépendance $(\star)$.\ebem
\bew D'après le Lemme \ref{L:CM=2ample}, il suffit de montrer que $T$ n'est pas $2$-ample au-dessus des $(T_i:i<n)$. Soient donc $\bar a$, $\bar b$ et $\bar c$ des uplets tels que~:
\begin{enumerate}
\item[(i)] $\bar a$ et $\bar b$ sont algébriquement clos~;
\item[(ii)] $\acl(\bar a,\bar b) \indi i_{\bar a} \sscl{\bar c,\bar a}$ pour tout $i<n$~;
\item[(iii)] $\bar c\ind_{\bar b}\bar a\bar b$.
\end{enumerate}
En ajoutant aux paramètres un modèle contenant $\bdd(\bar a)\cap\bdd(\bar b)$ et indépendant sur ces paramètres de $\bar a,\bar b,\bar c$,
on  peut supposer que $\bdd(\bar a)\cap\bdd(\bar b)=\bdd(\emptyset)$. Il suffira donc de montrer que $\bar c\ind\bar a$.
On supposera, grâce à la Remarque \ref{models}, que $\bar a$ est un modèle.

Posons $D=\sscl{\bar c,\bar a}\cap\acl(\bar c,\bar b)$. Nous allons montrer que $D\ind\bar a$ ce qui terminera la démonstration car $\bar c\in D$. Les conditions (ii) et (iii) impliquent que
$$\cb_i(D/\acl(\bar a,\bar b))\subseteq\bar a\cap\bar b=\acl(\emptyset)\quad\text{pour tout }i<n,$$
et donc $D\indi i\bar a$ pour tout $i<n$. D'après la caractérisation de l'indépendance $(\star)$, il reste à montrer que la sous-structure engendrée par $D\cup\bar a$ est autosuffisante, et  que ses points colorés sont les produits de ceux de $D$ et de $\bar a$ dans le cas modulaire.

La condition (ii) implique que $\sscl{\bar c,\bar a}\cap\acl(\bar a,\bar b)=\bar a$, donc 
$$\sscl{\bar c,\bar a}\cap\bar b\subseteq \bar a\cap \bar b=\acl(\emptyset).$$
Dans le cas dégénéré, le langage est purement relationnel et $D\cup\bar a$ est bien une sous-structure. Vérifions dans ce cas qu'elle est autosuffisante. L'indépendance (iii) implique par $(\star)$ que
$$\sscl{\acl(\bar c,\bar b),\acl(\bar a,\bar b)}=\acl(\bar c,\bar b)\cup\acl(\bar a,\bar b).$$
Ceci donne
$$\begin{array}{rcl}\sscl{\bar c,\bar a}\cap\sscl{\acl(\bar c,\bar b),\acl(\bar a,\bar b)}&=&
\sscl{\bar c,\bar a}\cap(\acl(\bar c,\bar b)\cup\acl(\bar a,\bar b))\\
&=&(\sscl{\bar c,\bar a}\cap\acl(\bar c,\bar b))\cup(\sscl{\bar c,\bar a}\cap\acl(\bar a,\bar b))\\
&=&D\cup\bar a.\end{array}$$
Comme l'intersection de deux ensembles autosuffisants l'est aussi, $D\cup\bar a$ est autosuffisante comme souhaité.

Considérons maintenant le cas modulaire. Notons que la structure $D\cdot\bar a$ engendrée par  $D \cup \bar a$ a pour domaine le produit des groupes $D$ et  $\bar a$ (y compris s'il y a des points de torsion). Rappelons que dans le cas coloré, seuls les points colorés peuvent avoir une prédimension relative négative. Par l'absurde prenons un uplet fini $\bar \gamma$ (de points colorés, dans le cas coloré) dans  $\sscl{D,\bar a}$ de longueur minimale satisfaisant l'une des deux conditions suivantes :
\begin{itemize}
\item[$(a)$] $\delta(\bar\gamma/D\cup \bar a)<0$,  ou
\item[$(b)$]$\bar\gamma$ est un point coloré de $D\cdot\bar a$ qui n'appartient pas au groupe engendré par les points colorés de $D \cup \bar a$.
\end{itemize}

La condition (iii) entraîne l'indépendance de $\acl(\bar a,\bar b)$ et $\acl(\bar c,\bar b)$ au-dessus de $\bar b$.  D'après la caractérisation de l'indépendance $(\star)$, le produit $\acl(\bar a,\bar b)\cdot \acl(\bar c,\bar b)$ des groupes $\acl(\bar a,\bar b)$ et $\acl(\bar c,\bar b)$ est autosuffisant et ses (éventuels) points colorés sont des produits de points colorés de $\acl(\bar a,\bar b) \cup \acl(\bar c,\bar b)$. Donc $\sscl{D,\bar a}\subseteq\acl(\bar a,\bar b)\cdot \acl(\bar c,\bar b)$ et $\bar\gamma=\bar\gamma_1\bar\gamma_2$ avec $\bar\gamma_1\in\acl(\bar a,\bar b)$, $\bar\gamma_2\in\acl(\bar c,\bar b)$ tel que, dans le cas coloré, $\bar\gamma_1$ et $\bar\gamma_2$ sont colorés. Nous allons montrer qu'on peut choisir $\bar\gamma_1 \in \bar a$ et $\bar\gamma_2 \in D$ ce qui contredira l'hypothèse sur $\bar \gamma$. 

Puisque $D\cup\bar\gamma\subseteq\sscl{D,\bar a}\subseteq\sscl{\bar c,\bar a}$, la condition (ii) entraîne pour tout $i<n$ l'indépendance 
\begin{equation}\label{eqn3}D,\bar\gamma\indi i_{\bar a}\acl(\bar a,\bar b).\end{equation}
Par la caractérisation $(\star)$, les ensembles $\acl(\bar c,\bar b)$ et $\acl(\bar a,\bar b)$ sont $i$-indépendants au-dessus de $\bar b$ pour tout $i<n$. Comme $D\cup\bar\gamma_2$ est contenu dans $\acl(\bar c,\bar b)$, on obtient, pour tout $i<n$, l'indépendance 
\begin{equation}\label{eqn4}D,\bar\gamma_2\indi i_{\bar b}\acl(\bar a,\bar b).\end{equation}
                                                                          
Pour un type $p$ stationnaire sur $A$ et $B\supseteq A$ on notera $p|_B$ l'unique extension non-déviante de $p$ sur $B$. Rappelons que deux types stationnaires $p$ et $q$ (éventuellement sur des  domaines différents) sont \emph{parallèles} (ce qu'on note $p \parallel q$) s'ils ont même extension non-déviante sur la réunion de leurs domaines. Le parallélisme est une relation d'équivalence. 

Pour chaque $i<n$, notons $p_i(X,\bar x,\bar a)=\tp_i(D,\bar\gamma/\bar a)$. Ces types sont stationnaires puisque $\bar a$ est un modèle. Pour $\bar\gamma'\in M^{|\bar\gamma|}$ et $\bar a'\models\tp(\bar a)$ on pose
$$\bar\gamma'\,p_i(X,\bar x,\bar a')=\{\p(X,\bar\gamma^{\prime-1}\bar x):\p(X,\bar x)\in p_i(X,\bar x,\bar a')|_{\bar a',\bar\gamma'}\},$$
le type qui correspond au translaté par $\bar\gamma'$ (sur les variables $\bar x$) de l'extension non-déviante de $p_i(X,\bar x,\bar a')$ à $\bar a',\bar\gamma'$. Autrement dit, 
$$\bar\gamma'\,p_i(X,\bar x,\bar a')= \tp_i(D_{\bar a'},\bar\gamma' \bar\gamma_{\bar a'} /\bar a', \bar\gamma')$$ pour $D_{\bar a'},\bar\gamma_{\bar a'}$ une réalisation de $p_i(X,\bar x,\bar a')$ tel que $D_{\bar a'},\bar\gamma_{\bar a'} \indi i_{\bar a'} \bar\gamma'$.

Par (\ref{eqn3}) et (\ref{eqn4}) remarquons que pour tout $i<n$, le type $\tp_i(D,\bar\gamma_2/ \bar b)$ est parallèle au translaté par $\bar\gamma_1^{-1}$ de $\tp_i(D,\bar \gamma/ \bar a)$. L'idée de la démonstration est de remplacer $\bar\gamma_1 \in \acl(\bar a,\bar b)$ par un uplet dans $\acl(\bar a)$. Pour cela on introduit la relation $E$ sur les réalisations de $\tp(\bar a)$ donnée par
$$\bar a'E\bar a''\quad\Leftrightarrow\quad\exists\,\bar\gamma'\in M^{|\bar\gamma|}\ \bigwedge_{i<n}\bar\gamma'\,p_i(X,\bar x,\bar a')\parallel p_i(X,\bar x,\bar a'').$$
Dans le cas coloré, on demande de plus que $\bar\gamma'$ soit coloré. 

Cette relation est type-définissable en utilisant la définissabilité des $i$-types.
Voyons que c'est une relation d'équivalence dont les classes sont  des  hyperimaginaires. La réflexivité est évidente. Montrons la transitivité : considérons $\bar a'$, $\bar a''$, $\bar a'''\models\tp(a)$ avec $\bar a'E\bar a''$ et $\bar a''E\bar a'''$, et choisissons $\bar\gamma'$,$\bar\gamma''$ (éventuellement colorés) avec
$$\bigwedge_{i<n}\bar\gamma'\,p_i(X,\bar x,\bar a')\parallel p_i(X,\bar x,\bar a'')
\quad\mbox{et}\quad
\bigwedge_{i<n}\bar\gamma''\,p_i(X,\bar x,\bar a'')\parallel p_i(X,\bar x,\bar a''').$$
Ceci signifie que pour tout $i<n$,
$$\begin{array}{rcl}\bar\gamma'\,p_i(X,\bar x,\bar a')|_{\bar a',\bar a'',\bar a''',\bar\gamma',\bar\gamma''}&=&
p_i(X,\bar x,\bar a'')|_{\bar a',\bar a'',\bar a''',\bar\gamma',\bar\gamma''}\\
\bar\gamma''\,p_i(X,\bar x,\bar a'')|_{\bar a',\bar a'',\bar a''',\bar\gamma',\bar\gamma''}
&=&p_i(X,\bar x,\bar a''')|_{\bar a',\bar a'',\bar a''',\bar\gamma',\bar\gamma''}.\end{array}$$
Donc
$$\begin{array}{rcl}\bar\gamma''\bar\gamma'\,p_i(X,\bar x,\bar a')|_{\bar a',\bar a'',\bar a''',\bar\gamma',\bar\gamma''} &=&
\bar\gamma''\,p_i(X,\bar x,\bar a'')|_{\bar a',\bar a'',\bar a''',\bar\gamma',\bar\gamma''}\\
&=&p_i(X,\bar x,\bar a''')|_{\bar a',\bar a'',\bar a''',\bar\gamma',\bar\gamma''}
\end{array}$$
et $\bar a'E\bar a'''$. La symétrie se montre de manière analogue.

Montrons maintenant que la classe de $\bar a$ modulo $E$ est bornée sur $\bar b$. On considère pour cela $\bar a',\bar\gamma'\models\tp(\bar a,\bar\gamma_1/\bar b)$ et on vérifie que $\bar a'$ est dans la classe d'équivalence de $\bar a$.
Fixons $i< n$. Puisque $\bar\gamma_1\in\acl(\bar a,\bar b)$ et  $D,\bar\gamma\indi i_{\bar a}\acl(\bar a,\bar b)$ d'après~(\ref{eqn3}), on~a 
\begin{equation}\label{eqn5}\bar\gamma_1^{-1}\,p_i(X,\bar x,\bar a)|_{\acl(\bar a,\bar b)}=\tp_i(D,\bar\gamma_1^{-1}\bar\gamma/\acl(\bar a,\bar b))=\tp_i(D,\bar\gamma_2/\acl(\bar a,\bar b)).\end{equation}
Par (\ref{eqn4})  ce type est l'unique extension non-déviante de $\tp_i(D,\bar\gamma_2/\bar b)$ à $\acl(\bar a,\bar b)$.
Comme  $\bar a',\bar\gamma'$ a même $i$-type que $\bar a,\bar\gamma_1$ sur $\bar b$,
$$\bar\gamma^{\prime-1}\,p_i(X,\bar x,\bar a')|_{\acl(\bar a',\bar b)}$$
est également l'unique extension non-déviante de $\tp_i(D,\bar\gamma_2/\bar b)$ à $\acl(\bar a',\bar b)$ et
$$\bar\gamma_1^{-1}\,p_i(X,\bar x,\bar a)|_{\acl(\bar a,\bar b)\cup \acl(\bar a',\bar b)}=
\bar\gamma^{\prime-1}\,p_i(X,\bar x,\bar a')|_{\acl(\bar a,\bar b)\cup \acl(\bar a',\bar b)}$$
est l'unique extension non-déviante de $\tp_i(D,\bar\gamma_2/\bar b)$ à $\acl(\bar a,\bar b)\cup \acl(\bar a',\bar b)$. On en conclut que $\bar aE\bar a'$ et donc que la classe de $\bar a$ modulo $E$ est bornée 
sur $\bar b$.

La classe de $\bar a$ modulo $E$ est contenue dans $\bdd(\bar a)\cap\bdd(\bar b)=\bdd(\emptyset)$.
Soit $N$ un modèle $\omega$-saturé tel que $\tp(N/\bar a,\bar b,\bar c)$ cohérite sur $\dcl(\emptyset)$ (qui est un modèle). Alors on trouve un représentant $\bar a_0$ dans $N$ de la classe $\bar aE$. Soit $\bar\gamma_0$ témoignant de l'équivalence de $\bar a_0$ et $\bar a$. Pour terminer il suffira de montrer que 
$$\bar\gamma_0\in\acl(\bar a,\bar a_0)\quad\mbox{et}\quad\bar\gamma_0^{-1}\bar\gamma_1\in\acl(\bar b,\bar a_0).$$
En effet, comme $\tp(\bar a_0/\bar a,\bar b)$ est finiment satisfaisable dans $\dcl(\emptyset)$ et $\bar\gamma_1\in\acl(\bar a,\bar b)$, on pourra trouver 
$\bar\gamma'_0\in\acl(\bar a)=\bar a$, éventuellement coloré, 
avec $\bar\gamma_0^{\prime-1}\bar\gamma_1\in\acl(\bar b)=\bar b$. Enfin
$$\bar\gamma_0^{\prime-1}\bar\gamma=\bar\gamma_0^{\prime-1}\bar\gamma_1\bar\gamma_2\in\sscl{D,\bar a}\cap\acl(\bar b,\bar c)\subseteq\sscl{\bar c,\bar a}\cap\acl(\bar b,\bar c)=D,$$
ce qui donnera une contradiction car $\bar \gamma$ sera alors le produit d'un élément (éventuellement coloré) de $\bar a$ et d'un élément (éventuellement coloré) de $D$.

Montrons d'abord que $\bar\gamma_0\in\acl(\bar a,\bar a_0)$. Rappelons que  $\bar\gamma_0$ vérifie  pour tout $i<n$,  
\begin{equation}\label{eqn6}\bar\gamma_0\,p_i(X,\bar x,\bar a_0)|_{\acl(\bar a,\bar a_0),\gamma_0} = p_i(X,\bar x,\bar a)|_{\acl(\bar a,\bar a_0),\gamma_0}.\end{equation}
Fixons une réalisation $D',\bar\gamma'$ de 
$\tp(D,\bar\gamma/\bar a)|_{\bar a,\bar a_0,\bar\gamma_0}$. Par l'hypothèse et (\ref{eqn6}), l'uplet $D',\bar\gamma_0^{-1}\bar\gamma'$ réalise $p_i(X,\bar x,\bar a_0)|_{\bar a,\bar a_0,\bar\gamma_0}$ pour tout $i<n$ et ainsi les uplets $D',\bar\gamma_0^{-1}\bar\gamma',\bar a_0$ et $D,\bar\gamma,\bar a$ ont mêmes $i$-types et même coloration. En particulier l'uplet $\bar\gamma_0^{-1}\bar\gamma'$ satisfait la même condition ($(a)$ ou $(b)$) sur $D',\bar a_0$ que celle satisfaite par $\bar \gamma$ sur $D ,\bar a$. On en déduit (à l'aide de la remarque \ref{sscl} pour le cas $(a)$) que $$\bar\gamma_0^{-1}\bar\gamma' \in \acl(D',\bar a_0)\quad\text{et donc}\quad\bar\gamma_0\in\acl(D',\bar a_0,\bar\gamma').$$
L'indépendance $\bar a_0,\bar\gamma_0\ind_{\bar a}D',\bar\gamma'$ implique alors que $\bar\gamma_0\in\acl(\bar a,\bar a_0)$.

Montrons enfin que $\bar\gamma_1^{-1}\bar\gamma_0 \in \acl(\bar b, \bar a_0)$. Comme $\bar a_0$ est indépendant de $\bar a, \bar b, D, \bar\gamma$, il suit par transitivité de l'indépendance et (\ref{eqn3}) que 
$$D,\bar\gamma\indi i_{\bar a}\acl(\bar a,\bar a_0,\bar b).$$
Cette indépendance, le fait que $\bar\gamma_0$ et $\bar\gamma_1$ sont dans $\acl(\bar a,\bar a_0, \bar b)$, et les égalités (\ref{eqn5}) et (\ref{eqn6}) permettent de vérifier que
$$\bar\gamma_1^{-1}\bar\gamma_0\,p_i(X,\bar x,\bar a_0)|_{\acl(\bar a,\bar a_0,\bar b)} =\bar\gamma_1^{-1}\,p_i(X,\bar x,\bar a)|_{\acl(\bar a,\bar a_0,\bar b)}=\tp_i(D,\bar\gamma_2/\bar b)|_{\acl(\bar a,\bar a_0,\bar b)}.$$
On considère alors une réalisation $D',\bar\gamma'$ de $\tp(D,\bar\gamma_2/\bar b)|_{\acl(\bar a,\bar a_0,\bar b)}$. D'après l'égalité précédente, l'uplet $D',\bar\gamma_0^{-1}\bar\gamma_1\bar\gamma'$ réalise $p_i(X,\bar x,\bar a_0)$ pour tout $i<n$, ce qui entraîne comme dans le paragraphe précédent que
$$\bar\gamma_0^{-1}\bar\gamma_1\bar\gamma'\in\acl(D',\bar a_0).$$
L'indépendance $\bar a_0,\bar\gamma_0^{-1}\bar\gamma_1\ind_{\bar b}D',\bar\gamma'$ implique alors que $\bar\gamma_0^{-1}\bar\gamma_1\in\acl(\bar b,\bar a_0)$.\ebew

\kor\label{C:corpscolores} Dans les corps colorés de rang de Morley fini tout groupe infini simple interprétable est linéaire. Aucun corps rouge n'admet de mauvais groupe interprétable. Si un mauvais groupe était interprétable dans un corps vert, il ne serait constitué que d'éléments semi-simples.\ekor
\bew Par le théorème précédent les corps colorés sont CM-triviaux au-dessus de la théorie des purs corps algébriquement clos (qui élimine les imaginaires). Considérons un groupe $G$ infini simple définissable dans une de ces théories. Par le Théorème \ref{T:homo}, il se plonge définissablement dans un groupe connexe algébrique $H$. Par le théorème de Chevalley \cite{Chev}, $H$ est une extension d'une variété abélienne (qui est commutative) par un groupe linéaire. Comme $G$ est simple, il doit se plonger dans la partie linéaire de $H$. En particulier, tout mauvais groupe interprétable est linéaire. Dans \cite{Polong} il est montré qu'un tel mauvais groupe ne peut exister qu'en caractéristique nulle et qu'il n'est constitué alors que d'éléments semi-simples (c.à.d.\ \emph{diagonalisables} en tant que matrices).\ebew

\bem
Dans un travail en cours \cite{BPW13}, les auteurs montrent que tout groupe simple définissable dans un corps coloré est définissablement isomorphe à un groupe algébrique. On vérifie également que tout groupe connexe définissable dans une fusion fortement minimale de deux théories fortement minimales avec élimination des imaginaires au-dessus de l'égalité est isogène à un produit de groupes définissables dans chacune des théories.
\ebem

\end{document}